\documentclass[11pt]{article}
\usepackage{amsmath,bm}            
\usepackage{array}              
\usepackage{latexsym}
\usepackage{epsfig}
\usepackage{makeidx}   
\usepackage{multirow}
\usepackage{comment}
\usepackage{booktabs} 
\usepackage{amssymb}
\usepackage{algorithm} 
\usepackage[noend]{algpseudocode}
\usepackage{graphicx}
\usepackage{xcolor} 
\usepackage{caption}

\graphicspath{ {./images/} } 
\usepackage{color,soul}
\allowdisplaybreaks

\usepackage{tikz,pgfgantt,setspace}
\newtheorem{proposition}{Proposition}
\newcommand*{\QEDB}{\hfill\ensuremath{\square}}%
\usepackage{url}
\begin{document}
	
	\title{A Stochastic Programming Approach for Chemotherapy Appointment Scheduling}
	\date\today

	\begin{titlepage}
		\vspace*{18ex} {\centerline{\bf{A STOCHASTIC PROGRAMMING APPROACH 
		}}} {\centerline{\bf{FOR CHEMOTHERAPY APPOINTMENT SCHEDULING}}}\vspace{0.3in}

		\begin{center} {\bf Nur Banu Demir$^1$, Serhat Gul$^2$, Melih \c{C}elik}$^{3}$ \\

			\vspace{1cm}
			$^1$Department of Industrial and Systems Engineering \\ Wayne State University, Detroit, Michigan, USA \\
			$^2$Department of Industrial Engineering, TED University, Ankara, Turkey \\ 
			$^3$School of Management, University of Bath, Bath, UK\\
			\textbf{Correspondence} \\
			Serhat Gul, Department of Industrial Engineering \\ TED University, 06420, \c{C}ankaya, Ankara, Turkey \\ 
			E-mail: serhat.gul@tedu.edu.tr 
		\end{center} \vspace{0.2in}
		
		\centerline{{July 2020}}
	\end{titlepage}
	
	\maketitle
	\begin{abstract}
		Chemotherapy appointment scheduling is a challenging problem due to the uncertainty in pre-medication and infusion durations. In this paper, we formulate a two-stage stochastic mixed integer programming model for the chemotherapy appointment scheduling problem under limited availability of nurses and infusion chairs. The objective is to minimize the expected weighted sum of nurse overtime, chair idle time, and patient waiting time. The computational burden to solve real-life instances of this problem to optimality is significantly high, even in the deterministic case. To overcome this burden, we incorporate valid bounds and symmetry breaking constraints. Progressive hedging algorithm is implemented in order to solve the improved formulation  heuristically. We enhance the algorithm through a penalty update method, cycle detection and variable fixing mechanisms, and a linear approximation of the objective function. Using numerical experiments based on real data from a major oncology hospital, we compare our solution approach with several scheduling heuristics from the relevant literature, generate managerial insights related to the impact of the number of nurses and chairs on appointment schedules, and estimate the value of stochastic solution to assess the significance of considering uncertainty.
		
	\end{abstract}
	
	\vspace{0.2in}
	
	\vspace{0.4in}
	\section{Introduction}
Recent statistical data on global health shows that cancer is one of the most prevalent causes of death, second only to heart disease (Ritchie and Roser, 2018)\nocite{cause}. The number of new cancer cases per year are estimated to increase from 17 million to 23.6 million by 2030 (Cancer Research UK, 2018; National Cancer Institute, 2018).\nocite{CancerUK} \nocite{NCI} In 2017, estimated expenditures for cancer treatment were \$147.3 billion in the US\nocite{NCI} (National Cancer Institute, 2018). Since many people are affected by this disease, it is crucial to have planning systems which decrease costs and increase patient satisfaction at the same time.
	
Chemotherapy is a frequently used method to cure cancer patients. Chemotherapy treatment is an exhausting process for patients, since side effects may be observed during or after the treatment. To reduce the intensity of the side effects, pre-medication drugs are injected to the patients before initializing the treatment. After the pre-medication process, the patients receive chemotherapy drugs. These drugs should be given to patients based on predetermined frequency and doses. 
	
Designing patient schedules is vitally important due to the limited capacity of oncology clinics and time restrictions for the chemotherapy treatment. Oncology clinics should provide schedules that consider the trade-off between provider and patient satisfaction. In the literature and in practice, chemotherapy schedules are generally created in two phases: (i) \textit{chemotherapy planning}, where patients are assigned to days according to their treatment frequency and (ii) \textit{chemotherapy scheduling}, where patients are sequenced and their appointment times are set to create a schedule for a given day. In this study, we concentrate on constructing daily schedules for chemotherapy patients. 
	
 Nurse overtime, chair idle time and patient waiting time are the three criteria that should be considered to evaluate the effectiveness of chemotherapy scheduling. Nurse overtime and chair idle time are undesirable for providers, since they increase operating costs of the clinic. Furthermore, working for more than the shift length may increase nurse dissatisfaction, and which in turn may lead to higher nurse turnover rates. Reducing patient waiting time is also important for clinics to improve patient satisfaction and service level. Since treatment durations are not actually known in advance of the treatment, rough estimates of these durations are generally used for scheduling. However, using such simple estimates may result in an undesirable amount of waiting time, idle time and overtime, particularly when the estimates are not very accurate.
	
The most complicating factor for chemotherapy patient appointment scheduling is the uncertainty in treatment durations. The treatment durations are unknown due to reasons including a change in the prescription list according to patient's health status before initiating the treatment, complications in patients, and early treatment termination due to the patient not tolerating the treatment (Castaing et al., 2016)\nocite{Castaing2016}. If the decision maker plans the daily schedule based on the longest possible treatment durations, patient waiting time can be reduced. However, this may lead to excessive nurse overtime and chair idle time. Similarly, schedules based on the shortest possible treatment durations may result in high waiting times in the clinic. Therefore, the decision maker should handle the uncertainty in treatment durations wisely to avoid long waiting time, chair idle time and nurse overtime. 
	
The limited availability of chairs and nurses in the system is another crucial factor in this problem. A patient should simultaneously seize a nurse and an infusion chair for the treatment.   The pre-medication process involves patients receiving pre-medication drugs, if required in the prescription. A nurse can conduct at most one pre-medication at a time. During chemotherapy infusion, the drugs that are used in chemotherapy treatment are injected to patients, in which case a single nurse proctors multiple patients. 
	
In this paper, we study the problem of sequencing patients and setting appointment times (i.e., scheduling appointment start times) for a chemotherapy unit considering the availability of a limited number of nurses and chairs under pre-medication and infusion duration uncertainties. Decisions in this problem include (1) sequencing patients of a daily appointment list, (2) setting appointment times, (3) assignment of patients to nurses, and (4) assignment of patients to chairs. We model the problem using a two-stage stochastic mixed integer programming (SMIP) formulation. We consider an objective function that minimizes the total expected penalty of patient waiting time, chair idle time and nurse overtime across a large set of scenarios sampled through pre-medication and infusion time distributions. We propose a modified progressive hedging algorithm (PHA) that utilizes the problem structure to find near-optimal chemotherapy schedules. In particular, we propose a penalty update method that considers convergence behavior of the primal and dual variables. The method includes a limit on the penalty parameter, whose value changes according to the iteration number. We also make use of cycle detection and variable fixing mechanisms, and linearize the objective function to improve solution times of scenario subproblems. Using instances based on data from a major oncology hospital, we compare the  performance of PHA with heuristics used in the relevant studies from the appointment scheduling literature. Our experiments provide insight into the issues related to the following questions:
	
	\begin{enumerate}
		\item What is the value of considering uncertainty in pre-medication and infusion durations when scheduling chemotherapy appointments?
		\item What is the potential benefit of using the PHA over commonly used heuristics from the relevant appointment scheduling literature?
		\item Which PHA routines and parameters must be more carefully designed to enhance the algorithm performance and solution quality?
		\item How does the number of nurses and chairs affect optimal schedules in terms of patient waiting time, chair idle time and nurse overtime?
	\end{enumerate}
	
	The organization of the remaining sections is as follows. In the next section, literature review and  the main contributions of the article are presented. In Section 3, the problem description and model formulation are given. In Section 4, the proposed progressive hedging algorithm and its implementation details are presented. Our computational experiments and their results are demonstrated in detail in Section 5. Finally, concluding remarks and further research directions are discussed in Section 6.
	
\section {Literature Review}
\label{section:2}
Our literature review mainly focuses on outpatient chemotherapy scheduling problems. The reader is referred to Lame et al. (2016) \nocite{Lame2016} for an extensive review on the subject. A detailed review of outpatient scheduling studies on primary care and specialty care can be found in Gupta and Denton (2008)\nocite{gupta2008appointment}, Cayirli and Veral (2003)\nocite{cayirli2003outpatient}, Ahmadi-Javid et al. (2017)\nocite{Ahmadi2017}. However, before we proceed with the review of chemotherapy scheduling literature, we discuss studies on general outpatient appointment scheduling. 

\subsection{General Outpatient Scheduling Problems}
Among the studies on outpatient scheduling, the ones that consider (i) patient sequencing (ii) appointment time decisions, or (iii) patient scheduling in a multiple resource setting are relevant to our study. As also pointed out in Ahmadi-Javid et al. (2017)\nocite{Ahmadi2017}, most of the articles assume a predetermined patient sequence due to the complexity of sequencing problem (Kong et al., 2019\nocite{Kongetal2019}; Begen and Queyranne, 2011\nocite{begenandqueyranne2011}; Denton and Gupta, 2003\nocite{dentonandgupta2003}). In the context of single-server scheduling, only a few articles consider sequencing and time setting decisions simultaneously (Berg et al., 2014\nocite{Bergetal2014}; Mak et al., 2014\nocite{Maketal2014}; Mancilla and Storer, 2012\nocite{Mancilla2012}; Denton et al., 2007\nocite{denton2007optimization}). Among these studies, Berg et al. (2014) \nocite{Bergetal2014}, Mancilla and Storer (2012)\nocite{Mancilla2012}, Denton et al. (2007)\nocite{denton2007optimization} formulated stochastic programming models. Note that the second-stage subproblems in those models are much simpler than those of the chemotherapy scheduling problem due to the consideration of two types of resources in the latter. 

In the literature of outpatient scheduling, the articles that consider multiple resources in the model are the most relevant ones to our study (Batun et al., 2011\nocite{Batunetal2011}; Perez et al., 2013\nocite{perezetal2013}, Riise et al., 2016\nocite{riiseetal2016}, Alvarez-Oh et al., 2018\nocite{alvarezohetal2018}, Klassen and Yoogalingam, 2019\nocite{KlassenYoogalingam2019}; Leeftink et al., 2019\nocite{Leeftinketal2019}; Hur et al., 2020\nocite{Huretal2020}).  Below, we provide an in-depth comparison between our study and the studies that formulate two-stage stochastic programming model.

Batun et al. (2011)\nocite{Batunetal2011} considered both surgeons and ORs in their study. However, the only assignment decision they made is the assignment of surgeries to ORs in their two-stage SP model.  In other words, they assumed that the surgery list of a surgeon is known. Furthermore, the surgery-to-OR assignment decision was made at the first stage of their SP model. This helped them model the second-stage problem using only continuous variables. Due to the resulting structure of the model, they were able to implement the L-shaped algorithm (Van Slyke and Wets, 1969\nocite{VanSlykeWets1969}) to solve their model. On the other hand, in our model we make assignments to two different types of resources including nurse and chairs. Furthermore, the assignment decisions are made at the second stage of the SP formulation. These variables make the structure of the second-stage problem of our model different from that of Batun et al. (2011)\nocite{Batunetal2011}. The types of variables at each stage of an SP formulation are critical in determining the candidate decomposition algorithms to be implemented. Due to the binary and continuous variables in the second stage, the L-shaped algorithm is not a viable solution method for our model. Hence, we implement a progressive hedging algorithm to obtain solutions. 

Perez et al. (2013)\nocite{perezetal2013} considered multiple equipment types and human resources while scheduling nuclear medicine appointments. They formulated both deterministic and stochastic formulations for the scheduling problem. They consider uncertainty in future patient arrival requests in their SP model but ignore any uncertainty in procedure durations. Note that the scheduling decisions in our SMIP formulation are made under uncertainty in treatment durations. Their SP formulation was designed to provide an appointment time for only one patient. To schedule multiple patients, they needed to solve the formulation multiple times within an algorithm. Their objective function includes the waiting time from the arrival time of the appointment request until the planned time of appointment. However, our model considers the waiting time from the time of patient arrival to clinic on the day of appointment until the actual treatment start time on the same day. Furthermore, we consider resource overtime as well. Perez et al. (2013)\nocite{perezetal2013} assigned patients to small time slots which allowed them to have only binary variables in both stages of the formulation. However, we do not use time slots in our formulation, and instead model appointment time decisions using continuous variables, resulting in a formulation with mixed binary and continuous variables in both stages. Having only binary variables in the first-stage problem would allow us to implement an integer L-shaped algorithm (Laporte and Louveaux, 1993\nocite{laporte1993integer}) or an evaluate-and-cut method proposed by Ahmed (2013)\nocite{Ahmed2013} to obtain optimal solutions. However, due to the aforementioned structure of our model, we implement a modified PHA. We do not use time slots, since modeling appointments using time slots would lead to unnecessary idle times.

Alvarez-Oh et al. (2018)\nocite{alvarezohetal2018} studied the allocation of patients to resources in a system with multiple service steps in the context of primary care clinic. However, they considered only patient-to-nurse assignments as the physician responsible for a patient is given in the problem. They also assigned patients to slots to determine appointment times, which allowed them to have only binary variables in the first stage of their SP formulation. Assignment of patients to multiple resources and the way appointment times are set constitute the main differences between our study and Alvarez-Oh et al. (2018)\nocite{alvarezohetal2018}. 

Hur et al. (2020)\nocite{Huretal2020} scheduled patients of a multidisciplinary outpatient clinic using SP formulation. They assigned patients to multiple providers to be seen during their visit at a back pain integrated practice unit. However, they assumed that each patient is seen by only one provider at a time. The appointment times were also set in terms of time blocks, which results in an SP formulation having only binary variables at the first stage. Note that in our case the patients are assigned to multiple resources simultaneously in the pre-medication step. Continuous appointment times make our SP model structure different from that of Hur et al. (2020)\nocite{Huretal2020} as well.

Leeftink et al. (2019)\nocite{Leeftinketal2019} also considered a multidisciplinary outpatient clinic to schedule appointments. However, their SP model is quite different from ours because they did not make individual patient assignments to resources. They instead created a blueprint schedule by determining the number of patients seen by each provider at each time period.

\subsection{Chemotherapy Scheduling Problems}
We do not restrict our review to studies that focus only on the treatment phase of chemotherapy (i.e. pre-medication and infusion). Some studies also consider the steps that need to be completed before the treatment starts, namely blood tests, oncologist evaluation, and drug preparation in the pharmacy. A patient's blood test result determines whether or not the patient is ready to go through the treatment. Once the blood test results are received from the laboratory, the treatment of the patient is either approved or postponed after the oncologist evaluates the results. In case the treatment is approved, a drug preparation order is provided to the pharmacist. The treatment can start when the chemotherapy drug is prepared.
	
Streams of literature most relevant to the work in this paper include the deterministic and stochastic versions of chemotherapy 
%planning and 
scheduling problems and their solution approaches, which are discussed in the remainder of this section.
\subsubsection{Deterministic Chemotherapy Scheduling Problems}
\label{sec:2.1}
Turkcan et al. (2012)\nocite{Turkcanetal2012} handled chemotherapy planning and scheduling problems in a hierarchical manner. In their first formulation, the aim is to minimize treatment delays while assigning the new patients to treatment days. After determining the daily patient lists, the authors focused on the daily patient scheduling problem by considering resource availabilities and acuity levels of the patients. Santibanez et al. (2012)\nocite{Santibanez2012} suggested that the planner should consider the preferences of patients, capacity of the pharmacy and laboratory, schedules of the oncologists, and overtime for nurses while assigning appointment times to patients. To determine the chemotherapy schedules, a timetable was prepared in Dobish et al. (2003)\nocite{dobish2003}, where it was assumed that laboratory tests and oncologist evaluation are completed on the previous day of the treatment. Pharmacist availability was an important concern in the scheduling process. Oncologist evaluation and pharmacy preparation stages were studied in addition to the treatment stage in the model developed in Sadki et al. (2011)\nocite{Sadki2011}. Nurses were always assumed to be ready for the treatment, so the scarce resources were the oncologists and chemotherapy chairs.

Heshmat et al. (2018)\nocite{heshmat2018solving} proposed a two-stage solution approach for the chemotherapy appointment scheduling problem. In the first stage, patients were grouped using clustering algorithms. In the latter stage, these clusters of patients were assigned to nurses, chairs, and time slots by solving a mathematical programming model. Huggins et al. (2014)\nocite{Huggins2014} developed a mathematical programming model to maximize chair utilization while considering the workloads of the pharmacists and nurses. They validated the model solutions with a simulation study. Two heuristics were developed to assign patients to the infusion chairs  by Sevinc at al. (2013)\nocite{Sevinc2013}. They decomposed the problem into two phases. They first determined the daily patient list for laboratory tests. The patients who obtain an oncologist approval according to their laboratory test results are then assigned to infusion chairs in the second phase. The makespan and weighted flow time were minimized in Heseraki et al. (2018)\nocite{hesaraki2018generating} to schedule appointments of the chemotherapy patients. The authors expressed treatment durations in terms of time slots by assuming that a nurse can proctor the pre-medication of only one patient at a given time. Other aspects considered in the study are the patient priorities and lunch-coffee breaks of the nurses.

Liang and Turkcan (2016)\nocite{Liangetal2016} discussed two different care delivery models for treatment of patients, namely \textit{functional} and \textit{primary care delivery models}. In the functional care delivery model, it is allowed to assign a different nurse to a patient at each treatment visit. On the other hand, a specific nurse is responsible for the patient's treatment in the primary care delivery model. The workload among nurses tends to show higher variability in the primary care delivery model, compared to the functional care delivery model. The authors constructed a mathematical programming model to minimize nurse overtime and total excess workload for the case of primary care delivery model. The acuity levels of patients are taken into account, and there is an upper bound on the assigned acuity level for each nurse in a time slot. Skill levels of the nurses are also included in the model. For the functional care delivery model, they aimed to minimize the total nurse overtime and patient waiting time within a multi-objective optimization model framework.
	
Our study differs from the studies summarized in this section since we consider uncertainty in the pre-medication and infusion durations.

\subsubsection{Stochastic Chemotherapy Scheduling Problems}
\label{sec:2.2}

Alvarado and Ntaimo (2018)\nocite{Alvarado2018} formulated three different mean-risk stochastic integer programming models for the chemotherapy appointment planning and scheduling problem. Acuity levels, treatment durations, and nurse availability were represented by stochastic parameters in their problem description. Their model was formulated to provide a solution for a single patient. In other words, the authors assumed that the appointments of earlier patients were already booked. They made their planning and scheduling decision for the current patient by considering the remaining time slots in the future days. The objective was to minimize the deviation from the target treatment start day for the new patient, patient waiting time, and nurse overtime. In our study, we consider planning the daily schedules of multiple patients without using time slots.

Hahn-Goldberg et al. (2014)\nocite{Hahn2014} considered the uncertainty due to the arrivals of appointment requests and last minute changes in a schedule. They minimized makespan in the dynamic scheduling problem using an online algorithm. However, they did not consider uncertainty in treatment durations. A discrete event simulation model was developed to observe the effect of different operational decisions on the patient waiting time, clinic overtime, and resource utilization in an outpatient oncology clinic in Liang et al. (2015)\nocite{Liang2015}. The stochastic elements in the simulation model included unpunctual arrivals of patients and all service durations. However, the uncertainty was not considered in the mathematical programming model used to create schedules.
	
Tanaka et al. (2011)\nocite{tanaka2011infusion} assigned patients to the infusion chairs using online bin packing heuristics. It was assumed that each infusion chair has a separate patient list for the treatments. The uncertainty in drug preparation and nursing durations before the infusion, such as taking vitals and assessment was taken into account. However, the study assumed deterministic pre-medication and infusion durations. Mandelbaum et al.\nocite{mandelbaumdata} (2016) considered uncertain treatment durations and unpunctual patients. They constructed a data-driven approach for the problem of patient appointment scheduling, assuming that the infusion chairs are the only servers of the system. Therefore, the impact of limited nurse availability on a schedule is not considered in this study.
	
The study most similar to that in this paper was conducted by Castaing et al. (2016)\nocite{Castaing2016}, where the authors constructed a two-stage stochastic integer program to determine patient appointment times in an outpatient chemotherapy clinic in order to minimize the expected patient waiting time and total time required for all treatments (i.e., clinic closure time). The authors assumed that an initial schedule was already created. Their model allowed to make small revisions on the existing schedule as they assumed the sequence of patients cannot be changed. They set appointment times for a given sequence of patients and considered multiple chairs, but assumed a single nurse in the clinic. The binary variables in the second stage made the model difficult to solve. Therefore, a two-phase heuristic which facilitates setting the values of those variables was implemented to solve the problem. Our work differs from Castaing et al. (2016)\nocite{Castaing2016} in a number of ways. First, we do not make the restricting assumption that the patient sequence is fixed. Hence, our model does more than the refinement of existing schedules; it can be used to create an initial schedule. Moreover, we consider that there are multiple nurses in the clinic and that a patient can be assigned to any of the available nurses for treatment. Therefore, our model is applicable to clinics operating based on a functional care delivery model. On the other hand, the model of Castaing et al. (2016)\nocite{Castaing2016} can be used by clinics functioning according to primary care delivery model. Furthermore, their model must be solved separately for each nurse, as it considers a single nurse. 
	
\section{The Chemotherapy Appointment Scheduling Problem}
\label{chp:b3}
In this section, we define the Chemotherapy Appointment Scheduling Problem based on our observations of the chemotherapy operations at the Hacettepe University Oncology Hospital in Ankara, Turkey, which is one of the prominent oncology centers in the country. We also describe the proposed two-stage stochastic mixed integer programming formulation in detail.

\subsection{Chemotherapy Operations in the Hacettepe Outpatient Chemotherapy Unit}
The Hacettepe Outpatient Chemotherapy Unit consists of 32 chemotherapy chairs, 28 of which are used for treatments. The remaining four chairs are reserved for supportive care. The number of patients receiving treatment varies between 60 and 85 on a given day. On average, 10 nurses, including a head nurse, provide service to patients each day. During chemotherapy planning (which is outside the scope of this paper), patients are assigned to each day based on the day of their previous visit, estimated duration of their treatment, and the expected treatment duration of patients assigned to each upcoming day. A daily shift starts at 8:00 and ends at 17:00 (unless overtime is needed), with a lunch break between 12:00-13:00. The daily schedule of patients is arranged by the head nurse, who assigns appointment times, chairs and nurses to the patients.
	
Our observations at the chemotherapy unit aimed to capture the operations conducted until and during the chemotherapy treatment. The blood tests and oncologist evaluations are carried out on a different day before the treatment. The results of these tests are examined by the oncologist to decide whether the patient has convenient health status to receive a chemotherapy treatment. When the oncologist approves the treatment protocol, she updates or confirms the dosages of the drugs. The drugs are then prepared in the pharmacy lab before the patient arrives at the hospital for her treatment. During each treatment visit, a patient receives two types of medications: (i) pre-medication drugs, which are injected prior to the chemotherapy infusion to help prevent side effects; and (ii) chemotherapy infusion drugs, which are used for cancer treatment.
	
There are certain time slots that the patients are assigned to according to their estimated treatment duration, which is predicted by the chemotherapy unit according to the types and dosages of the drugs injected. During a day, the unit makes use of four time slots: 8:00-10:30, 10:30-12:00, 13:00-15:30 and 15:30-17:00. If the estimated treatment duration of a patient is less than the length of a slot, the patient is given an appointment time so that the expected end of treatment is within the same slot. Otherwise, the patient can be assigned to more than one slot. If the estimated treatment duration of a patient is longer than those of other patients, the head nurse prefers to assign this patient to the beginning of the workday in order to prevent excessive overtime. In the rare case that the estimated duration of the treatment exceeds 270 minutes, the appointment is generally split into two sessions, one in the morning and one in the afternoon. Given the uncertainties in treatment times, the three main objectives in designing the schedule are to ensure that (i) each patient can actually start their treatment after experiencing minimal or no wait, (ii) all chairs are highly utilized, and (iii) each nurse spends minimal or no time above regular shift duration.
	
\begin{figure}
		\includegraphics[width=\linewidth]{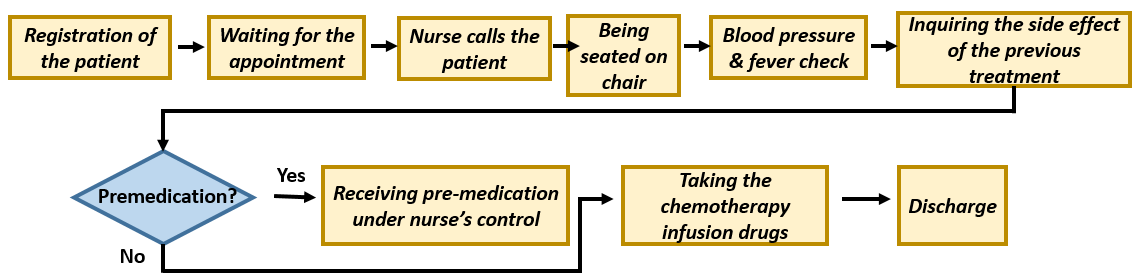}
		\caption{Patient flow chart for treatment process}
		\label{fig:Picture1}
\end{figure}

The treatment of a patient consists of a sequence of events, as illustrated in Figure \ref{fig:Picture1}.  A patient that arrives at the unit is immediately registered and becomes available for treatment at the appointment time. If a nurse and/or chair is not available at the time, the patient has to wait for the treatment. A nurse calls the patient when an infusion chair and the nurse both become available. After the patient is seated in a chemotherapy chair, the nurse measures the fever and blood pressure, and establishes a vascular access to start the pre-medication process. Aside from rare cases where the patient may not need pre-medication drugs, the infusion process starts after the pre-medication process ends. When the infusion process is over, the patient is discharged, and the nurse calls the next patient to a chair that has become available.
In general, each nurse may be responsible for up to four patients at a given time. The Hacettepe Outpatient Chemotherapy Unit applies a modified form of the \textit{functional care delivery}, where in line with this procedure each patient can be assigned to different nurses in subsequent visits. However, nurse-patient assignments in the clinic also consider the type of cancer and expected duration of treatment, which may limit the pool of nurses that can be assigned to the patient. The balance of both criteria is aimed to be satisfied to create equity between nurses.
\subsection {A Two-Stage Stochastic Mixed-Integer Programming Model for the Chemotherapy Appointment Scheduling Problem}
We formulate a two-stage stochastic mixed-integer programming (SMIP) model for the Chemotherapy Appointment Scheduling Problem, motivated by the operations at the Hacettepe Outpatient Chemotherapy Unit. In particular, we study the problem of sequencing patients and setting appointment times for a chemotherapy unit by considering the availability of nurses and chairs under uncertainty on the pre-medication and infusion durations. We sequence patients of a daily appointment list, set appointment times, and assign patients to nurses and chairs in our model.  We consider an objective function that minimizes a weighted combination of patient waiting time, chair idle time and nurse overtime.

In the remainder of this section, we provide our assumptions and two-stage stochastic programming model and propose a number of valid inequalities and bounds to strengthen the model. 
\subsubsection{Two-Stage Stochastic Mixed-Integer Programming Formulation}	
In this section, we discuss our two-stage SMIP model for the chemotherapy appointment scheduling problem.	The following assumptions are made in the model:
	\begin{itemize}
		\item We take the decisions in the chemotherapy planning phase (appointment of patients to days) as given and focus on the chemotherapy scheduling decisions for each day.
		\item Among the events shown in Figure \ref{fig:Picture1}, we focus on the three main events in between the patient appointment time and treatment completion time. These events include patient waiting, pre-medication, and infusion. 
	    \item We do not consider time slots (as applied in the Outpatient Chemotherapy Unit) while determining the patient appointment times. Although the use of slots facilitates the construction of schedules manually, having to fit the treatment of patients into pre-determined slots may be too restrictive, as it may lead to long idle times within some slots. To overcome this, our model schedules patients in any minute within each half-day shift.
		\item Patients can be treated by any available nurse, as is the case for \textit{functional care delivery}. Nurses are also assumed to have identical skills.
		\item A nurse can perform only one patient's pre-medication process at any time.
		\item A nurse can proctor the infusion process of multiple patients while conducting the pre-medication process of a patient.
		\item Patients become ready exactly on the appointment time for treatment.
	   	\item The sequence of patients is not allowed to change after they arrive at the chemotherapy unit. 
	\end{itemize}
At the first stage of the model, patients are sequenced and their appointment times are set.  These constitute the main decisions in the model. Then, pre-medication and infusion durations are realized. At the second stage, patients are assigned to nurses and chairs. The remaining decision variables at the second stage help determine patient waiting time, chair idle time, discharge time, and nurse overtime. Overtime is calculated separately for each nurse that works for more than the planned shift length. Chairs become idle due to two reasons. First, treatment of a patient may finish earlier than the appointment time of the subsequent patient. Second, longer than expected durations of pre-medications may prevent nurses from initiating the treatment of a waiting patient. Waiting occurs when the treatment starts later than the appointment time of the patient, which may happen when pre-medications or infusions for the previous patients may last longer than expected. The unavailability of nurses or chairs may also cause waiting. 

\begin{table}
\centering
\caption{Notation used throughout the two-stage SMIP model}
\label{tab:notation}
\small
\begin{tabular}{cl}
\multicolumn{2}{c}{\textbf{Index sets}} \\
\hline
$I$ & Patients \\ 
$C$ & Chairs	\\ 
$\Omega$ & Scenarios \\
$N$ & Nurses \\
\hline
& \\
\multicolumn{2}{c}{\textbf{Parameters}} \\
\hline
$s^{\omega}_{i}$ & Pre-medication time of patient $i \in I$ in scenario $\omega \in \Omega$ \\ 
$t^{\omega}_{i}$ & Infusion time of patient $i \in I$ in scenario $\omega \in \Omega$ \\
$p^{\omega}$ & Probability of scenario $\omega \in \Omega$ \\
$H$ & Shift duration of nurses \\
$L$ & Overtime limit of nurses \\
$\lambda_{j}$ & Trade-off parameters in the objective function in the interval $[0,1],\, j \in \{1,2,3\}$ \\
$M$ & A large value\\ 
\hline
& \\
\multicolumn{2}{c}{\textbf{First-stage decision variables}} \\
\hline
$b_{ij}$ &  $= \begin{cases} 1, & \text{if patient $i \in I$ precedes patient $j \in I$ in the daily appointment list } \\ 0, & \text{otherwise} \end{cases}$ \\
$a_i$ & Appointment time of patient $i \in I$ \\ 
\hline
& \\
\multicolumn{2}{c}{\textbf{Second-stage decision variables}} \\
\hline
$x^{\omega}_{in}$  & $= \begin{cases} 1, & \text{if patient $i \in I$ is assigned to nurse $n \in N$ in scenario $\omega \in \Omega$} \\ 0, & \text{otherwise} \end{cases}$ \\
$y^{\omega}_{ic}$  & $= \begin{cases} 1,  & \text{if patient $i \in I$ is assigned to chair $c \in C$ in scenario $\omega \in \Omega$} \\0, & \text{otherwise} \end{cases}$\\
$w^\omega_{i}$ & Waiting time of patient $i \in I$ in scenario $\omega \in \Omega$ \\ 
$d^\omega_{i}$ & Discharge time of patient $i \in I$ in scenario $\omega \in \Omega$\\
$O^\omega_{n}$ & Overtime of nurse $n \in N$ in scenario $\omega \in \Omega$\\
$I^\omega_{c}$ & Idle time of chair $c \in C$ in scenario $\omega \in \Omega$\\
$T^\omega_{c}$ & Auxiliary variable used to calculate idle time of chair $c \in C$ in scenario $\omega \in \Omega$\\
\hline
\end{tabular}
\end{table}

We assume a finite set of scenarios representing uncertainty in pre-medication and infusion durations. Given these scenarios, we formulate the following two-stage SMIP model for our problem based on sets, parameters, first and second-stage decision variables summarized in Table \ref{tab:notation}:
	{\begin{flalign}
		min \quad & \mathcal{Q}(\textbf{a,b})  & \label{obj}\\
		& b_{ij} + b_{ji} = 1   \quad \quad & \forall i, j \in I, j>i \label{prec}\\
		&a_j \geq a_i  -M (1-b_{ij}) \quad \quad & \forall i, j \in I, j \neq i \label{start}\\
		& b_{ij} \in \{0, 1\}  \quad \quad & \forall i, j \in I, j \neq i \label{bsign}\\
		& a_i :integer  \quad \quad & \forall i \in I \label{asign}
		\end{flalign}
		\begin{flushleft}
			$\text{where} $ \\
			$\quad \quad \mathcal{Q}(\textbf{a,b}) =  E_{\xi}[\text{Q}(\textbf{a,b},\xi(\omega))] $ \\
			$\text{is the expected recourse function, and} $ \\
			$\quad \quad \text{Q}(\textbf{a,b},\xi(\omega)) = \text{min}\bigg\{\lambda_{1} \displaystyle\sum_{i\in I}  w^\omega_{i}+\lambda_{2} \displaystyle\sum_{n\in N} O^\omega_{n} + \lambda_{3} \displaystyle\sum_{c\in C} I^\omega_{c}\bigg\}$ \\	
		\end{flushleft}
		\begin{flalign}
		&\displaystyle\sum_{n\in N} x^{\omega}_{in} =1 \quad & \forall i \in I \label{x}\\
		&\displaystyle\sum_{c\in C} y^{\omega}_{ic} =1 \quad & \forall i\in I \label{y}\\
		&a_i + w^\omega_{i} + s^\omega_{i} + t^\omega_{i} = d^\omega_{i}  \quad & \forall i\in I \label{timeeq}\\
		&a_j + w^\omega_{j} \geq a_i + w^\omega_{i} + s^\omega_{i}-M (3-b_{ij}- x^{\omega}_{in}- x^{\omega}_{jn}) & \forall i, j\in I,j\neq i,\forall n\in N \label{nurseprec}\\ 
		&a_j + w^\omega_{j} \geq d^w_{i}-M (3-b_{ij}-y^{\omega}_{ic}-y^{\omega}_{jc}) \quad &\forall i, j \in I,  j \neq i, \forall c \in C \label{yprec} \\ 
		&a_j + w^\omega_{j} \geq a_i + w^\omega_{i} -M (1-b_{ij}) \quad & \forall i, j \in I,  j \neq i \label{bprec}\\ 
		&O^\omega_{n} \geq d^\omega_{i} -H- M (1-x^{\omega}_{in}) \quad & \forall i \in I,\forall n \in N \label{overtime}\\
		&O^\omega_{n} \leq L \quad & \forall n \in N \label{imp3}\\
		&T^\omega_{c} \geq H \quad & \forall c \in C \label{idle1}\\
		&T^\omega_{c} \geq d^\omega_{i}-M(1-y^{\omega}_{ic}) \quad &\forall i \in I, \forall c \in C \label{idle2}\\
		&I^\omega_{c} \geq T^\omega_{c}-\displaystyle\sum_{i\in I} (s^\omega_{i}+t^\omega_{i}) y^{\omega}_{ic}  \quad & \forall c \in C \label{idle3}\\
		&x^{\omega}_{in} \in \{0,1\}  \quad & \forall i \in I,\forall n\in N  \label{xsign} \\
		&y^{\omega}_{ic} \in \{0,1\}  \quad & \forall i\in I,\forall c\in C \label{ysign} \\
		&d^\omega_{i},  w^\omega_{i} \geq 0  \quad & \forall i \in I \label{dwsign}\\
		&O^\omega_{n} \geq 0  \quad & \forall n \in N \label{Osign}\\
		&T^\omega_{c},  I^\omega_{c} \geq 0  \quad & \forall c \in C \label{Isign}
		\end{flalign}
	}
%\pagebreak\\
The objective function in (\ref{obj}) only includes the expected second-stage function, since there is no contribution from the first stage to the objective function. The expected second-stage function aims to minimize the expected weighted sum of total patient waiting time, nurse overtime, and chair idle time. First-stage constraints are given by \eqref{prec}-\eqref{asign}. Patient precedence relations are determined using constraints \eqref{prec}. If patient $i \in I$ is scheduled before the patient $j \in I$, the corresponding $b_{ij}$ value should be equal to 1. Constraints \eqref{start} establish the relationship between the binary precedence variables and appointment times, in that if patient $i \in I$ precedes patient $j \in I$ in the list, then the appointment time of patient $i \in I$ should be no later than that of patient $j \in I$. Binary and integrality restrictions on the first-stage variables are represented in constraints \eqref{bsign} and \eqref{asign}, respectively.
	
%In the second stage, nurse and chair assignments are made based on the realized treatment durations. 
For each scenario representing a set of chemotherapy durations, a subproblem is solved to obtain second-stage decisions. The constraints for the second-stage scenario subproblems are expressed by \eqref{x}-\eqref{Isign}. Constraints \eqref{x} and \eqref{y} enforce every patient to be assigned to exactly one nurse and one chair, respectively. Constraints \eqref{timeeq} calculate the patient discharge time, which is equal to the sum of the patient appointment time, waiting time, and durations of pre-medication and infusion. Constraints \eqref{nurseprec} ensure that if patient $i \in I$ and $j \in I$ are assigned to the same nurse, and patient $i \in I$ is scheduled before $j \in I$, then the treatment start time of patient $j \in I$ should start after the end of pre-medication of patient $i \in I$. Similar relation also holds for chair-precedence relationship, which is demonstrated in constraints \eqref{yprec}. If two patients are assigned to the same chair, the treatment of the latter may begin only after the discharge time of the former one. 
By constraints \eqref{bprec}, a patient's treatment start time should be earlier than those of the patients that are scheduled later. 
Overtime of a nurse should be either zero or the difference between the discharge time of the last discharged patient assigned to this nurse and the shift length, as determined by constraints \eqref{overtime}. In oncology units, it is desired to have limited overtime to enhance nurses' working conditions and preserve equity among nurses. A predetermined bound is included to ensure this in the model, which also tightens the feasible region of the formulation. The associated constraint is provided in equation \eqref{imp3}. Constraints \eqref{idle1}-\eqref{idle3} are used to calculate idle time of the chairs in the clinic. The idle time of each chair is calculated considering the shift length and discharge time of the lastly scheduled patient who is treated in the associated chair. The maximum of the discharge time of the last patient and shift length is used to measure the idle time of the associated chair. The remaining constraints are the sign and binary restrictions on the second-stage decision variables.

Given that nurses and chairs are identical, assigning patients to them in the first stage may lead to inefficient schedules. In practice, assignment of each patient to chairs and nurses is generally made dynamically after arrival, well before the uncertainty on the pre-medication and infusion durations is fully realized. Although treating these second-stage decisions in our model may seem in contradiction to practice, we show that once the patient sequencing is fixed in the first stage, these decisions become straightforward. More formally, we use the following proposition, the proof of which is given in Appendix \ref{sec:proof}.

\begin{proposition}
	Given the sequencing of patients from the first stage, assigning arriving patients to the first available chair and nurse produces the optimal schedule for the second stage of the Chemotherapy Appointment Scheduling Problem.
\end{proposition}

Using this proposition, we are able to consider chair and nurse assignment decisions in the second stage without violating the dynamic nature of the optimal assignment policy. Note that a particular assumption, which is also represented by constraints  \eqref{bprec}, has a critical role in the proof of the proposition. Constraints \eqref{bprec} enforce that the appointment sequence in the first stage is the same as the treatment sequence in the second stage. This ensures that the optimal policy to assign patients to chairs and nurses is non-idling. Hence, there is no need to make such assignments dynamically when patients arrive.    

The rule dictated by the above discussed assumption is valid from a practical point of view. As  all patients have the same urgency levels and are assumed to arrive on time, the managers would not ask an arriving patient to wait and let his/her successor be treated first when a nurse and a chair are available in the unit. This approach would dramatically reduce the satisfaction of the waiting patient. Therefore, it may be impractical to reshuffle patient sequence after the patient arrivals.  

We note here that although it is desired, we do not model a limit on the number of patients to be observed by a nurse at a given time explicitly. This potential issue is implicitly addressed by the model with the inclusion of the first two components of the objective function, as assigning too many patients to a single nurse at a given time may excessively increase the patient waiting time and nurse overtime. Hence, these two objective components would drive the solution to a balanced assignment of patients to nurses. In case an explicit limit $B$ on the number of patients assigned to a nurse at a time is desired, the following constraints can be added to the two-stage SMIP model:
\begin{align}
& u_{ink}^\omega \geq x_{in}^\omega + \beta_{ink}^\omega + \gamma_{ink}^\omega -2 & \forall i \in I,\, n \in N,\, k \in K\\
& M\beta_{ink}^\omega \geq k-a_i-w_i^\omega & \forall i \in I,\, n \in N,\, k \in K \\
& M\gamma_{ink}^\omega \geq a_i+w_i^\omega+s_i^\omega+t_i^\omega-k & \forall i \in I,\, n \in N,\, k \in K \\
& \sum_{i \in I} u_{ink}^\omega \leq B & \forall n \in N,\, k \in K\\
& u_{ink}^\omega,\, \beta_{ink}^\omega,\, \gamma_{ink}^\omega \in \{0,1\} & \forall i \in I,\, n \in N,\, k \in K,
\end{align}
where $K$ represents the set of discrete time units (e.g., minutes), $\beta_{ink}^\omega$ takes a value of 1 if patient $i \in I$ has started treatment with nurse $n \in N$ before $k \in K$, $\gamma_{ink}^\omega$ is 1 if patient $i \in I$ has not finished treatment with the same nurse before $k \in K$, $u_{ink}^\omega$ represents whether patient $i \in I$ is being treated by nurse $n \in N$ at time $k \in K$, and $M$ is a large number. It should be noted here that such a modeling approach would require discretization of the time periods (e.g., in minutes), and would substantially add to the complexity of the model. Furthermore, Proposition 1 should be modified to assign arriving patient to the first available nurse observing the treatment of less than $B$ patients.

\subsubsection{Improvements in the Formulation}
\label{3.c}
Based on the preliminary runs, it is observed that solving even deterministic instances of the problem necessitates significantly long solution times. Consequently, symmetry-breaking constraints and valid bounds are added in order to strengthen the formulation and decrease the computational burden. 

The first set of constraints relates to breaking symmetry with respect to the overtime of nurses. Since nurses are assumed to be identical, the schedules of the nurses are interchangeable. By convention, we enforce the overtime of nurse with index 1 to be no less than those of the other nurses. In this way, overtime may be assigned to nurses in non-increasing order of indices. The mathematical representation of this constraint is given by:
\begin{align}
&O^\omega_{n} \geq O^\omega_{n+1} \quad & \forall n \leq |N|-1, \forall \omega \in \Omega \label{imp1}
\end{align}
Note that balancing nurse workload is not among the main objectives in our model. However, the upper limit introduced by constraint \eqref{imp3} defines an acceptable bound and hence prevents excessive imbalance. %Therefore, adding constraint \eqref{imp1} into the model results in an implementable solution. Furthermore, the preferred nurses for overtime work can be changed each day, since the model proposes a schedule for a single day.

Next, a lower bound on nurse overtime is defined using the structure of the problem. If one can equally distribute the summation of treatment durations to chairs, the resulting average total treatment time per chair provides the smallest possible value for the maximum of patient discharge times. Therefore, the calculated average value is a lower bound for the sum of nurse shift time and overtime of nurse 1, since the overtime of nurse 1 is at least as high as all others according to equation \eqref{imp1}. The mathematical representation of this constraint is shown given by:
\begin{align}
&O^\omega_{1} +H \geq \frac{{\sum_{i \in I}} (s^\omega_{i}+t^\omega_{i})}{|C|}  \quad &\forall \omega \in \Omega \label{imp4} 
\end{align}
As proposed in Mancilla and Storer (2012\nocite{Mancilla2012}), the appointment and waiting times of the first patient in the schedule can be assigned to zero to reduce the number of decision variables. According to the results of preliminary experiments, adding such time assignments does not have a significant effect to reduce the computation times, and thus these are not included in the formulation. 

Equations \eqref{imp1}-\eqref{imp4} are added to the main model and used in the computational experiments.
\section{Solution Methodology}
\label{chp:b4}

SMIP models are in general computationally demanding mainly due to the large amount of variables and constraints depending on the number of scenarios (Birge, 2011\nocite{birge2011introduction}). Decomposition based solution methods are generally utilized to solve such models in the literature.  The well-known stage decomposition algorithm proposed by Slyke and Wets (1969)\nocite{van1969shaped}, \textit{L-shaped method}, cannot be applied on stochastic programs in case there are binary variables in the second stage as in our SMIP formulation. To handle those cases, Laporte and Louveaux (1993)\nocite{laporte1993integer} extended the method and proposed the \textit{integer L-shaped algorithm}, which can be implemented when only binary variables exist in the first stage. Other solution methods aim to overcome the difficulty associated with having binary variables in the second stage through approaches based on value function and set convexification (Sen, 2005;\nocite{sen2005c3} Yuan, 2009\nocite{yuan2009enhanced}). These disjunctive decomposition-based branch-and-cut algorithms are known to be very effective, but they can also be applied on models with only binary variables in the first stage. On the other hand, in our formulation, the first stage includes both binary and general integer variables.

Scenario decomposition-based approaches are also frequently used to solve SMIP models. The Progressive Hedging Algorithm (PHA) has been increasingly preferred among them (Gul et al., 2015\nocite{Gul2015}; Hvattum et al., 2009\nocite{Hvattum2009}; Goncalves et al., 2012\nocite{Goncalves2012}). The PHA, proposed by Rockafellar and Wets (1991)\nocite{Rockafellar1991}, is based on an augmented Lagrangian relaxation technique. The algorithm decomposes the problem into single-scenario subproblems, and aims to obtain an \textit{admissible} solution that is feasible for all scenarios by aggregating the scenario solutions at each iteration. It requires reformulation of the model to attain a structure appropriate for scenario decomposition. It then relaxes the constraints that enforce the first-stage decisions to be the same for every scenario in the reformulation and penalizes the violation of these constraints by introducing a Lagrangian term and quadratic penalty term in the objective function. The algorithm is shown to converge to the optimal solution for convex models. Since our model includes general integer and binary decision variables, it is not convex. Therefore, PHA is used as a heuristic approach for the Chemotherapy Appointment Scheduling Problem. Note that there is significant evidence in the literature showing that the PHA is an effective heuristic for SMIP models (Crainic et al., 2011\nocite{Crainic2011}; Gul et al., 2015\nocite{Gul2015}; Watson and Woodruff, 2011\nocite{Watson2011}).

In the remainder of this section, general steps of the PHA are given, followed by a literature review on different applications of the PHA. We make use of a number of approaches from the literature to modify the PHA, which is discussed in the last part of this section.
\subsection{The Progressive Hedging Algorithm}
\label{sec:4.1}
We first reformulate our SMIP model to  facilitate scenario decomposition. In the reformulation for PHA, which we call as SMIP-R (see Appendix \ref{sec:DeterministicEquivalentModel} for the model), \textit{non-anticipativity constraints} are explicitly formulated. To formulate these constraints, a copy of a first-stage variable must be created for each scenario. The constraints then enforce the first-stage decision variables to be equal to each other across all scenarios to prevent anticipation of the future. Since the appointment time values ($a_{i}$) and the precedence relationship variables ($b_{ij}$) are dependent on each other, there is no need to include non-anticipativity constraints for the latter. In SMIP-R, non-anticipativity constraints for the appointment time values are given by:
\begin{align}
&a^\omega_{i} = a_{i} \quad \forall i \in I, \forall \omega \in \Omega,
\end{align}
where the $a_{i}$ values represent the \textit{consensus variable} for appointment time of patient $i \in I$. Prior to the PHA implementation, these non-anticipativity constraints are relaxed, and their violation is penalized using two terms in the objective function. In the first term, the amount of violation for each patient is multiplied by Lagrangian multipliers denoted by $\mu^{\omega}_{i} \mbox{ } \forall i \in I,\forall \omega \in \Omega$. The second term is a quadratic component that takes the squared sum of these violations, and multiplies it by $\dfrac{\rho_{1}}{2}$, where $\rho$ denotes the penalty parameter. The revised objective function is given by:
\begin{comment}
\begin{align}
 &\lambda_{1} \displaystyle\sum_{i\in I} \displaystyle\sum_{\omega \in \Omega}p^\omega  w^\omega_{i} + \lambda_{2}\displaystyle\sum_{n\in N}\displaystyle\sum_{\omega \in \Omega} p^\omega  O^\omega_{n} +
 \lambda_{3}\displaystyle\sum_{c\in C}\displaystyle\sum_{\omega \in \Omega} p^\omega  I^\omega_{c} + \displaystyle\sum_{i\in I} \displaystyle\sum_{\omega \in \Omega}\mu^\omega_{i}(a^\omega_{i}-a_{i}) \nonumber & \\ 
 &+\displaystyle\sum_{i\in I} \displaystyle\sum_{n\in N}\displaystyle\sum_{\omega \in \Omega}\theta^\omega_{in}(x^\omega_{in}-x_{in})+\displaystyle\sum_{i\in I} \displaystyle\sum_{c\in C} \displaystyle\sum_{\omega\in\Omega}\gamma^\omega_{ic}(y^\omega_{ic}-y_{ic})+\dfrac{\rho_{1}}{2}\displaystyle\sum_{i\in I}\displaystyle\sum_{\omega \in \Omega} ||a^\omega_{i}-a_{i} ||^2 \nonumber& \\
 &+\dfrac{\rho_{2}}{2}\displaystyle\sum_{i\in I}\displaystyle\sum_{n\in N}\displaystyle\sum_{\omega \in \Omega} ||x^\omega_{in}-x_{in} ||^2 + \dfrac{\rho_{3}}{2}\displaystyle\sum_{i\in I}\displaystyle\sum_{c\in C}\displaystyle\sum_{\omega \in \Omega} ||y^\omega_{ic}-y_{ic} ||^2\label{zmu} 
\end{align}
\end{comment}
\begin{align}
&\lambda_{1} \displaystyle\sum_{i\in I} \displaystyle\sum_{\omega \in \Omega}p^\omega  w^\omega_{i} + \lambda_{2}\displaystyle\sum_{n\in N}\displaystyle\sum_{\omega \in \Omega} p^\omega  O^\omega_{n} +
\lambda_{3}\displaystyle\sum_{c\in C}\displaystyle\sum_{\omega \in \Omega} p^\omega  I^\omega_{c} + \displaystyle\sum_{i\in I} \displaystyle\sum_{\omega \in \Omega}\mu^\omega_{i}(a^\omega_{i}-a_{i}) \nonumber & \\ 
&+\dfrac{\rho_{1}}{2}\displaystyle\sum_{i\in I}\displaystyle\sum_{\omega \in \Omega} ||a^\omega_{i}-a_{i} ||^2  &\label{zmu} 
\end{align}

To obtain a separable formulation, variables $a_i$ in \eqref{zmu} are replaced by $\hat{a_i}$, a consensus parameter that is equal to the weighted average of appointment times over all scenarios, where the weights are equal to the probabilities of scenarios.
\begin{align}
&\hat{a_{i}} = \displaystyle\sum_{\omega \in \Omega}p^\omega a^\omega_{i} \quad \forall i \in I 
\end{align} 

The resulting structure of the objective function and constraint set allow the SMIP-R to be decomposed into multiple scenario subproblems, which are solved independently at each iteration of the algorithm.

The general steps of the PHA, shown also in Algorithm \ref{alg:PH}, are as follows: In step 1, algorithm parameters are initialized. Penalty parameter of the quadratic term, denoted by $\rho$, is set as a positive constant value of  $\rho^0$, and the Lagrangian multipliers ${\mu^{\omega (v)}_{i}}$  %\mbox{ } \forall i \in I,\forall \omega \in \Omega$ 
are initialized as 0 in the first iteration. In step 2, each scenario subproblem is solved. Note that quadratic penalty and Lagrangian  terms are ignored in the first iteration of this step. In step 3, the appointment time values obtained for every patient and scenario are used to calculate the consensus parameters $\hat{a_{i}}$ for every patient. %$i \in I$
In step 4, the penalty parameter ${\rho}$ is updated by multiplying it with a pre-determined positive constant ${\alpha}$. By using the value of the penalty parameter, Lagrangian multipliers are updated considering the distance between scenario solutions and the consensus parameter value. At the end of each iteration, the termination criterion is checked. If all first-stage solutions are identical, the algorithm stops. Otherwise, the algorithm returns to step 2 and solves scenario subproblems with the updated values of the multipliers and consensus parameter. 

\subsection{Previous Studies on the Progressive Hedging Algorithm}
\label{sec:4.2}
Hvattum et al. (2009)\nocite{Hvattum2009}, Gul et al. (2015)\nocite{Gul2015},  Gade et al. (2016)\nocite{Gade2016}, Watson et al. (2011)\nocite{Watson2011}, Cheung et al. (2015)\nocite{Cheung2015}, Mulvey and Vladimirou (1991)\nocite{Mulvey1991},  Goncalves et al. (2012)\nocite{Goncalves2012},  Atakan et al. (2018)\nocite{atakan2018progressive}, and Helgason et al. (1991)\nocite{helgason1991approximate} used PHA in application areas such as inventory-routing, surgery scheduling, unit commitment, hydrothermal systems, and fisheries management. A number of modifications to improve convergence behavior and computation times of the PHA were suggested by Watson and Woodruff (2011)\nocite{Watson2011}, who classified the changes made in the algorithm in four categories: computing effective $\rho$ values, accelerating convergence, termination criteria and detecting cyclic behaviour. Cheung et al. (2015)\nocite{Cheung2015} solved scenario subproblems approximately in the early iterations of the PHA, since obtaining exact solutions from the scenario subproblems is not a necessity for the algorithm. PHA was extended with some innovations such as dynamic multiple penalty parameters and heuristic intermediate solutions in Hvattum et al. (2009)\nocite{Hvattum2009}, where the penalty parameter is increased when the distance between the scenario solutions and the consensus parameter increases. When the consensus parameters at consecutive iterations moves farther away from each other, the penalty parameter is reduced.  Mulvey and Vladimirou (1991)\nocite{Mulvey1991} also applied a dynamic penalty parameter update method. The technique can prohibit stalling and ill-conditioning solution structures, which occurs when the algorithm gets stuck at suboptimal solutions.

\begin{algorithm}
	\caption{Progressive Hedging Algorithm}\label{alg:PH}
	\begin{algorithmic}[1]
		\State \textbf{Step 1: Initialization:}
		\State  Let $v$=1
		\State\indent ${\rho_{1}}^{(v)} = {\rho_{1}}^0$
		\State\indent $\mu^{\omega (v)}_{i} = 0 \quad \forall i \in I, \forall \omega \in \Omega$
		\State \textbf{Step 2: Solving scenario-subproblems:}
		\State\textbf{if} $v$ $=$ 1; 
		\State $\mbox {For each }  \omega \in \Omega; \mbox{ignore Lagrangian and quadratic penalty terms and solve  } $
		\State $\mbox{ scenario subproblems to obtain:}$
		\State\indent ${a^\omega_{i}}^{(v)} \quad\forall i \in I, \forall \omega \in \Omega $
		\State\textbf{end if}
		\State\textbf{else}
		\State $\mbox {For each }  \omega \in \Omega; \mbox{ solve scenario subproblems to obtain:} $ 
		\State\indent ${a^\omega_{i}}^{(v)} \quad\forall i \in I, \forall \omega \in \Omega $
		\State\textbf{end else}
		\State \textbf {Step 3: Calculate the consensus parameters:} 
		\State \indent $\hat{a_{i}} = \displaystyle\sum_{\omega \in \Omega}p^\omega a^{\omega (v)}_{i} \quad \forall i \in I$
		\State \textbf {Step 4: Update penalty parameters:}
		\State\textbf{if}  $v$ $>$ 1;  
		\State\indent ${\rho}^{(v+1)} = \alpha {\rho}^{(v)} $
		\State\textbf{end if}
		\State \textbf{Step 5: Update Lagrange Multipliers:}
		\State\indent ${\mu^\omega_{i}}^{(v+1)} = {\mu^\omega_{i}}^{(v)} +{\rho}^{(v)}({a^\omega_{i}}^{(v)}-\hat{a_{i}}) \quad \forall i \in I, \forall \omega \in \Omega$
		\State \textbf{Step 6: Termination:}
		\State \textbf{if} all first stage solutions are the same, then stop.
		\State\indent$ {a^{\omega (v)}_{i}} = \hat{a_{i}}  \quad \forall i \in I, \forall \omega \in \Omega$ 
		\State\textbf{end if}
		\State \textbf{else}
		\State\indent Set $v$= $v$+1; 
		\State \textbf{end else}
		\State \textbf{return} to Step 2
	\end{algorithmic}
\end{algorithm}

\subsection{Modifications on the Progressive Hedging Algorithm}
When the PHA is applied to the Chemotherapy Appointment Scheduling Problem, even single-scenario subproblems require substantial computational times, which deem realistically-sized instances of the problem impossible to solve in reasonable time. A number of issues contribute to this issue. The first one is the quadratic term used in the objective function of the subproblems, whereas the second one is the choice of penalty parameter update method. We investigate other improvement ideas on the PHA by detecting and favoring solutions found in the majority of the scenario subproblems, detecting cycles, and varying the optimality gap parameter of CPLEX to obtain approximate scenario subproblem solutions to decrease solution times. We also propose a lower bounding scheme for nurse overtime by benefiting from the structure of the PHA. In the following sections, these modifications on the PHA are presented. 
\subsubsection{Handling the Quadratic Term in the Objective Function}
The scenario subproblems solved at each iteration of the PHA includes a quadratic term in the objective function, which makes them difficult to solve. This term can be expanded as:
\begin{flalign}
\dfrac{\rho}{2}\displaystyle\sum_{i\in I}\displaystyle\sum_{w \in \Omega} ||a^{\omega}_{i}-\hat{a_{i}} ||^2 = \dfrac{\rho}{2}\displaystyle\sum_{i\in I}\displaystyle\sum_{\omega \in \Omega}(a^{\omega}_{i})^2 -\rho \displaystyle\sum_{i\in I}\displaystyle\sum_{\omega \in \Omega} (a^{\omega}_{i}\hat{a_{i}}) + |\Omega|\dfrac{\rho}{2}\displaystyle\sum_{i\in I}\hat{a_{i}}^2 \label{split}
\end{flalign}
Since $\hat{a_{i}}$ is a parameter, the only quadratic part is the first term on the right-hand side of \eqref{split}. Our aim is to use a method to handle this term in a fast manner without sacrificing significantly from the solution quality. For this end, we experiment with three different approaches: (1) nonlinear optimization using commercial solvers, (2) second-order cone programming (SOCP), and (3) linearization of the quadratic term.

For the first approach, Nonlinear Optimization package of CPLEX 12.8 is used. The quadratic term is included in the objective function without making any change in the algorithm. In the second approach, additional variables are defined to convert the quadratic term into an appropriate format for second-order cone programming. For the third approach, two alternatives are available in the PHA literature. The quadratic term is linearized using a piecewise linear function in Watson et al. (2012)\nocite{Watson2012}. Alternatively, Helseth (2016)\nocite{Helseth2016} dynamically approximated the quadratic term by benefiting from tangent line approximation. In particular, linear cuts that approximate the quadratic term are added to scenario subproblems to obtain a linear objective function. 

Based on our preliminary experiments, we use the third approach, since it outperforms the other two in terms of the computational times significantly, without substantial sacrifice in terms of solution quality. Due to the use of linearized penalty components in the objective function, we call our algorithm as the \textit{linearized progressive hedging algorithm} (LPHA). The linearization method that we use, which is based on the approach in Helseth (2016)\nocite{Helseth2016}, is explained next. 

A general form of Taylor's approximation, $f(m) \approx f(n) + f'(n) (m-n)$, can be used to derive a cut that provides a lower bound for the nonlinear term. If the cut is generated at different operating points (n) iteratively, a close approximation of $(a^{\omega}_{i})^2$ is obtained after some iterations. 
\begin{comment}
\begin{align}
&f(m) 	\approx f(n) + f'(n) (m-n)& \label{taylor}
\end{align}   
\end{comment}

In our case, at iteration $v$, the value of $(a^{\omega}_{i})^2$ is approximated based on the operating point $(a^{\omega (v-1)}_{i})$ which represents the appointment time set for patient $i$ in scenario $\omega$ at the previous iteration. Note that an alternative operating point would be selected as $\hat{a_{i}}$ calculated in iteration $(v-1)$. However, we choose the former alternative based on the results of preliminary experiments. After replacing $(a^{\omega}_{i})^2$ by a new variable $g^{\omega}_{i}$, the proposed linear cuts are then formulated as follows:
\begin{align}
&(g^{\omega}_{i}) \geq (a^{\omega (v-1)}_{i})^2+2(a^{\omega (v-1)}_{i})(a^{\omega}_{i}-a^{\omega (v-1)}_{i}) & \label{cut}
\end{align}
The cuts in \eqref{cut} are added to scenario subproblems at each iteration until it is noticed that the value of $g^{\omega}_{i}$ does not change from one iteration to another.   
\subsubsection{Choice of the Penalty Parameter Update Method}
In the PHA, the penalty parameter $\rho$ is updated by multiplying it with constant $\alpha$ in every consecutive iteration. When the multiplier $\alpha$ is set to 1, the penalty parameter becomes stable in the algorithm. Previous studies have shown that low values of $\rho$ tend to improve the quality of the solutions at the expense of longer computational times. On the other hand, high values of the penalty parameter can lead to oscillations in the solutions and the algorithm generally converges faster to a low-quality solution. Therefore, it is crucial to find a balance between solution time and quality using an appropriate level of the penalty parameter over iterations. Designing a well-performing penalty update method is challenging due to this trade-off. Besides, there is no straightforward way to select the initial value of the penalty parameter; it is mostly dependent on the objective function coefficients and scales of the decision variable values in the problem. 

One way to resolve the trade-off between solution quality and computational time is to apply a dynamic penalty parameter scheme by updating the value of the parameter in subsequent iterations. In this study, we use a method inspired by the approach discussed in Hvattum and Lokketangen (2009)\nocite{Hvattum2009}. In this method, two parameters are made use of by exploiting the relationship between the primal and dual. The first parameter, denoted by $\Delta_p$ measures the square of the difference between consensus parameters for the appointment time of each patient at consecutive iterations. If $\Delta_p$ is increasing, this means that consensus parameter is changing at a greater rate. This also implies that the current value of consensus parameter is not promising, and therefore it is risky to enforce convergence immediately. The second parameter, denoted by $\Delta_d$, compares the convergence behaviour of the first-stage variable to the consensus parameter at consecutive iterations. If $\Delta_d$ is increasing, this indicates that the first-stage variables in different scenario subproblems are farther away from the consensus parameter. 

In our implementation of the penalty parameter update method, when $\Delta_d$ increases in consecutive iterations, penalty parameter $\rho$ is multiplied by a constant $\alpha$ that is greater than 1. By this way, the algorithm is more inclined to converge faster. If the difference between $\Delta_d$ in consecutive iterations is not positive, the behaviour of $\Delta_p$ is observed. If $\Delta_p$ is increasing, the penalty parameter is reduced by a factor of $\frac{1}{\alpha}$. If both consecutive $\Delta_p$ and $\Delta_d$ differences are negative, the current penalty parameter is preserved. 

Based on preliminary results, we observe that penalty parameter can quickly climb to very large values as a result of multiplying with $\alpha$. In terms of convergence, higher penalty parameter values lead the algorithm to prematurely converge to suboptimal solutions. To avoid this issue, the value of penalty parameter is bounded above by a pre-determined parameter $\rho^u$, so that the algorithm can converge to a higher quality solution. 

It is also observed in our preliminary experiments that keeping the penalty parameter small is beneficial in terms of solution quality. However, if the predefined bound for the penalty parameter is kept smaller than needed during the iterations, this may result in excessive computational times. To avoid this, we preserve a small bound at the earlier iterations to obtain high quality solutions, and then we increase the limit in the later iterations to achieve fast convergence. The details of the penalty update method is provided in Algorithm \ref{alg:PU4}. In the algorithm, $\rho_1^u$ and $\rho_2^u$ refer to two different penalty parameter limits, where $\rho_1^u < \rho_2^u$. Furthermore, $iterlimit$ corresponds to the iteration value at which the upper bound is changed from $\rho_1^u$ to $\rho_2^u$.
\algtext*{Indent}
\begin{algorithm}
	\caption{Penalty Update Method with Limit}\label{alg:PU4}
	%\vspace{-12pt}
	%\doublespacing
	\begin{algorithmic}[1]
		\State ${\Delta_p}^{(v)} = {\displaystyle\sum_{i\in P}}(\hat{a_{i}}^{(v)}-\hat{a_{i}}^{(v-1)})^2 $
		\State ${\Delta_d}^{(v)} ={\displaystyle\sum_{i\in P}}{\displaystyle\sum_{\omega\in \Omega}}({a^{\omega (v)}_{i}}-\hat{a_{i}}^{(v)})^2 $
		\State\textbf{if } $v\leq iterlimit$
		\State \indent $\rho^u = \rho_1^u$
		\State\textbf{else }
		\State \indent $\rho^u = \rho_2^u$
		\State\textbf{if}  ${\Delta_d}^{(v)} -{\Delta_d}^{(v-1)} > 0 \quad \& \quad {\rho}^{(v)} < \rho^u;$
		\State\indent	${\rho}^{(v+1)} = \alpha{\rho}^{(v)}$
		\State\textbf{elseif}  ${\Delta_d}^{(v)} -{\Delta_d}^{(v-1)} > 0 $
		\State\indent	${\rho}^{(v+1)} = \rho^u $
		\State\textbf{elseif}  ${\Delta_p}^{(v)} -{\Delta_p}^{(v-1)} > 0 $
		\State\indent	${\rho}^{(v+1)} = \frac{1}{\alpha}{\rho}^{(v)}$
		\State\textbf{elseif} ${\rho}^{(v)} \leq \rho^u;$
		\State\indent ${\rho}^{(v+1)} = {\rho}^{(v)}$
		\State\textbf{else}	${\rho}^{(v+1)} = \rho^u$
	\end{algorithmic}
\end{algorithm}
\subsubsection{Cycling and Majority Value Detection}
Cycles are frequently observed within a typical implementation of the PHA. This behaviour may prevent convergence. We benefit from the idea of cycling prevention suggested by Watson and Woodruff (2011)\nocite{Watson2011} to overcome this. The appointment time values may appear identical over several consecutive PHA iterations. The direct consequence of this is that the consensus parameter values remain constant over the iterations. On the other hand, Lagrangian multipliers change in every iteration, since we update their values using Step 5 of Algorithm 1.

While updating the Lagrangian multipliers, the differences between appointment time values and consensus parameter are penalized by multiplying it with a penalty parameter. Due to this reason, Lagrangian multipliers vary over iterations even if the appointment time and consensus parameter values stay the same. Hence, any cycling behaviour of the algorithm can be detected by examining the behaviour of the Lagrangian multipliers.

Watson and Woodruff (2011)\nocite{Watson2011} used a hashing scheme to detect cycles. In our case, we would like to detect cycles in integer appointment times for each patient. Therefore, we have to check the behavior of the Lagrangian multipliers associated with each patient. In their approach, Watson and Woodruff (2011)\nocite{Watson2011} first generated hash weights for each scenario ($z_{\omega}$). In our algorithm, hash weights are the random numbers between 0 and 1. These hash weights are then multiplied with the Lagrangian multipliers ($\mu^{\omega}_{i}$). Hash value for each patient is calculated in every iteration as:
\begin{flalign}
h_{i} = \sum_{\omega \in \Omega} z_{\omega}\mu^{\omega (v)}_{i} \quad \forall i \in I \label{hash}
\end{flalign}
Note that $h_{i}=0$ at iteration $v=1$, as $\mu^{\omega (1)}_{i}=0$. If $z_{\omega}$ values are set equal to the probability of occurrence of each scenario ($p^\omega$), then $h_{i}$ would be zero for $v>1$ due to the reason briefly explained next. Since ${\mu^\omega_{i}}^{(v+1)} = {\rho}^{(v)}({a^\omega_{i}}^{(v)}-\sum_{\omega \in \Omega}p^\omega a^{\omega (v)}_{i})$ for all $i \in I, \omega \in \Omega$ at $v=1$, $\sum_{\omega \in \Omega} p^{\omega}\mu^{\omega (2)}_{i} = 0$. By induction, it can also be verified that $\sum_{\omega \in \Omega}p^{\omega}\mu^{\omega (v)}_{i} = 0$ for all v. Therefore, $z_{\omega}$ is not set as $p^{\omega}$ in our experiments.

To detect cycling, we compare values of $h_{i}$ in consecutive iterations. Since the continuous values of hash weights and Lagrangian multipliers may result in non-integer hash values, we compare the hash values of patients in consecutive iterations by using a small threshold value. If the difference between the hash values are smaller than the threshold value over three consecutive iterations, it is considered as the indicator of a cycle. The cycle then needs to be broken to maintain algorithm convergence. 

After a cycle is detected, the appointment time for patient $i$ can be fixed to a certain value. There might be different approaches for fixing the appointment time, which include the minimum and maximum appointment times observed among all scenario subproblem solutions for that patient in the previous iteration. Based on our preliminary experiments, we fix the appointment time to the one that is most repeatedly observed among subproblem solutions, which we call the \textit{majority value}. We use this approach after the first 50 iterations of the PHA, since there is high variation in the appointment times of a patient across different scenarios in the beginning iterations. The coupled implementation of cycle detection and variable fixing methods ensures the convergence of the PHA to a local optimum in finitely many steps.

A common observation in our preliminary experiments is that for many patients, the appointment times for a majority of the scenarios converge quickly, whereas those of the remaining scenarios fail to converge for many iterations. To overcome this, we fix the appointment time of a patient if there is convergence in 80\% of the scenarios. Although this may reduce the solution quality, it significantly reduces computational times. 

\subsubsection{Varying CPLEX optimality gap value in the subproblems}
\label{gap}
In the earlier iterations of the PHA, computation times to solve the subproblems are longer, since the algorithm uses lower values of penalty parameter $\rho$. However, the PHA does not need optimal solutions of the subproblems for convergence. We may heuristically solve the subproblems to decrease solution times. For this purpose, we test the idea proposed by Watson and Woodruff (2011)\nocite{Watson2011}, where the CPLEX optimality gap value for the subproblems is varied according to the convergence behaviour. In particular, it is monotonically decreased by checking the following condition:
\begin{flalign}
gap_{v+1} =\mbox{min \Bigg\{} {\sum_{i \in I, \hat{a_{i}}^{(v)}\neq 0} \sum_{\omega \in \Omega} \dfrac{|a^{\omega (v)}_{i} -\hat{a_{i}}^{(v)}|}{\hat{a_{i}}^{(v)}} ,\mbox{  } gap_{v}} \mbox{\Bigg\}}
\end{flalign}
The gap value is decreased if the appointment times of the patients move closer to the consensus parameter. Otherwise, we keep the gap value equal to that of the previous iteration. The initial gap value is determined according to the results of the preliminary experiments.

Using approximate solutions within the PHA may either improve solution quality by making greater amount of iterations or reduce it because of the low-quality solutions particularly at the earlier iterations. The details of this procedure and results are discussed in Section \ref{chp:b5}.

\subsubsection{Lower Bounds on Total Nurse Overtime and Chair Idle Time}
In the first iteration of the PHA, the scenario subproblems are solved without adding any Lagrangian or quadratic terms. Since the durations are known in each subproblem, waiting times for all patients would be reduced to zero without causing an increase in the total overtime or idle time values. In other words, the trade-off between waiting time and the sum of other measures disappears as the non-anticipativity constraint is totally omitted under this setting. Then, nurse overtime and chair idle time remain as the challenging terms in the subproblem objective function. Note that these two measures do not conflict with each other. Therefore, the resulting sum of nurse overtime and chair idle time values in a scenario subproblem solution in the first iteration constitutes a lower bound for the associated scenario subproblem in the subsequent iterations.

\section{Computational Experiments}
\label{chp:b5}
In this section, we discuss the computational experiments for the Chemotherapy Appointment Scheduling Problem. We test the methods using data gathered from the Outpatient Chemotherapy Unit at Hacettepe Oncology Hospital from November 2017 to March 2018. The data set includes realized pre-medication and infusion times as well as estimated treatment durations for all patients.

In our experiments, Microsoft Visual C++ 2019 is used for the implementation of the algorithm using CPLEX 12.9 Concert Technology. The computations are performed with Intel (R) Core (TM) i7-9750H CPU @2.60 GHz and 16GB RAM. We conduct each experiment on a problem instance set that consists of 10 instances. We generate those problem instances by sampling pre-medication and infusion durations in our data set.

In Section \ref{sec:IG}, the descriptive statistics on the data set is provided and the problem instance generation approach is discussed. The method proposed to improve scenario subproblem solution times is tested in Section \ref{testCPLEX}. 
%The sensitivity of the linearized progressive hedging algorithm (LPHA) to the algorithm parameters is investigated in Section \ref{sec:CR}.
 The performance of the LPHA is assessed in Section 5.3, whereas sensitivity analysis on model parameters is conducted to generate managerial insights in Section 5.4. Finally, the value of stochastic solution (VSS) is estimated in Section 5.5. 

\subsection{Data Description and Problem Instance Generation}
\label{sec:IG}
Our data set consists of the estimated and actual pre-medication and treatment durations of 204 patients we have collected data on over 11 different days. The total treatment duration for these patients ranges between 16 and 240 minutes. In particular, pre-medication duration varies between 0 and 36 minutes with the average of 15.3 minutes. Infusion duration changes between 16 and 217 minutes, and the average is 112.5 minutes. Furthermore, a 95\% confidence interval for the means of pre-medication and infusion durations are [14.61, 16.09] and [104.56, 120.35] minutes, respectively. 

\begin{figure}
\centering
	\includegraphics[height=7cm]{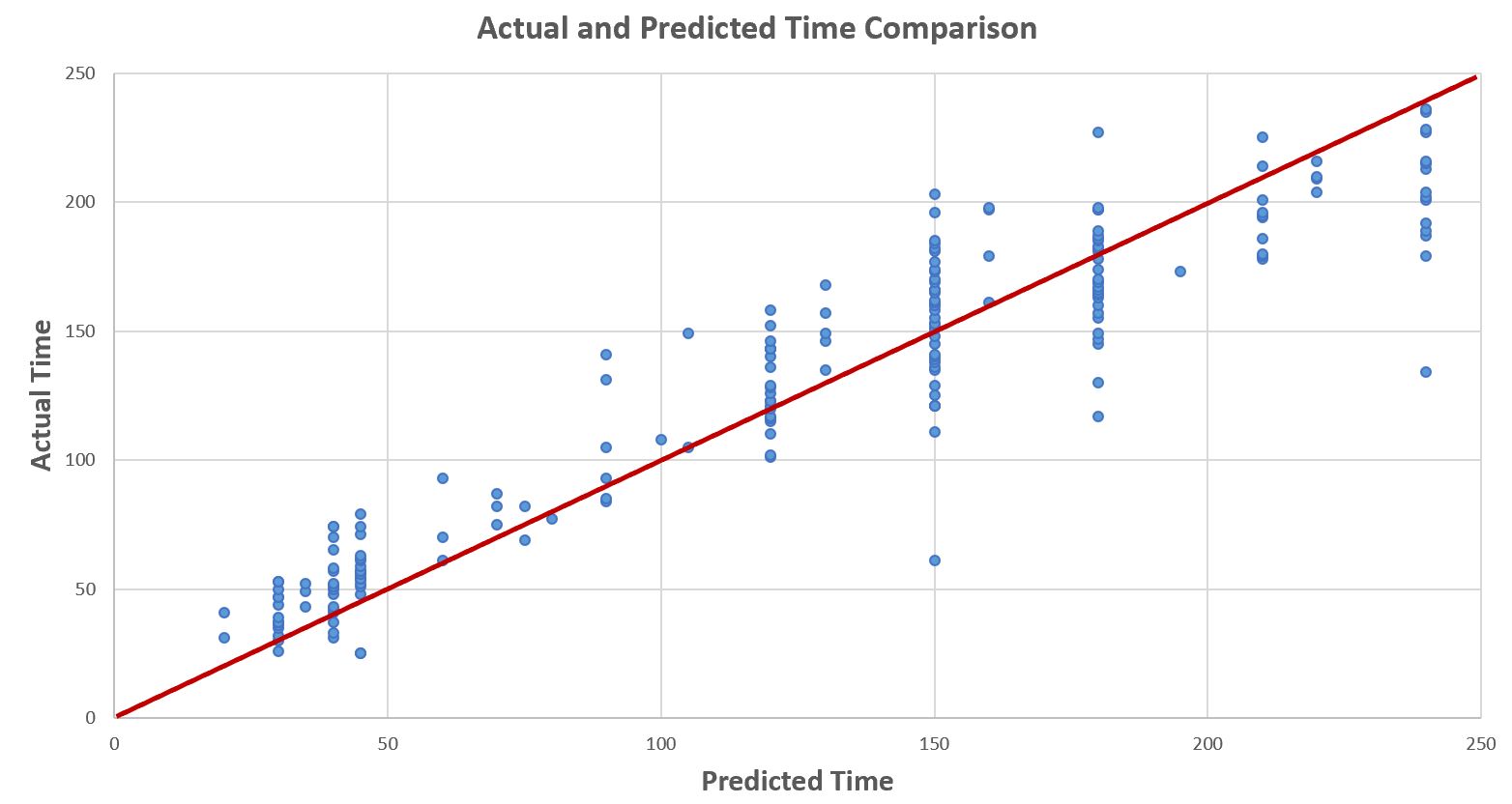}
	\caption{Comparison of actual and predicted treatment times}
	%\caption{Patient flow chart for treatment process}
	\label{fig:Picact}
\end{figure}

We next compare the realized treatment durations with the estimated durations. The data set indicates that 59\% of the patients were allocated less time than needed, whereas the remaining 41\% were allocated more than the necessary amount. Furthermore, the distribution of patient groups (by estimated duration) varies significantly among different half-shifts. To avoid sampling bias due to this reason, we generate representative half-shifts in our instances, rather than replicating any of the half-shifts. In Figure \ref{fig:Picact}, the comparison of actual and predicted treatment times can be observed on a scatter plot. In horizontal and vertical axes, the predicted and actual treatment times are provided, respectively. The value of mean absolute percent error (MAPE) is 16.9\%, which shows a significant variability in observed treatment time values compared to predicted values. Furthermore, Figure \ref{fig:Picact} also indicates an underestimation of the actual times of treatment times shorter than 150 minutes, and overestimation of those longer than three hours. Therefore, the manager of the chemotherapy clinic should consider a scheduling model where the variability of treatment times are taken into account.

In the Hacettepe Outpatient Chemotherapy Unit, the schedules are designed for half-day periods. To mimic the current practice followed in the unit, we also design half-day schedules. We generate each problem instance by sampling pre-medication and infusion durations from the data set. As can also be seen in Figure \ref{fig:Picact}, there are 21 discrete groups based on the estimated treatment duration of the patient. Since each half-shift consists of 7–8 patients and some of the 21 groups have very limited number of observations, we cluster these into four larger classes with cut-offs of 45, 100, 150 and 240 minutes of estimated treatment time. For each class, the actual pre-medication and infusion durations are analyzed. The intervals for the realized pre-medication and infusion durations for the defined classes as well as the probability of observing the defined classes in the data set are provided in Table \ref{durations}. 
\begin{table}
	\centering
	\caption{Treatment duration intervals (in minutes) and their frequencies in the data set}
	\begin{tabular}{|c|c|c|c|c|}
		\hline
		\textbf{Patient} & \textbf{Estimated} & \textbf{Associated} & \textbf{Actual pre-medication} & \textbf{Actual infusion} \\
		\textbf{Class} & \textbf{interval} & \textbf{probability} & \textbf{interval} & \textbf{interval} \\
		\hline
1	&	[20, 45] & 26.96\% & [0, 14] & [16, 44] \\ \hline
2	&	(45, 100] & 7.85\% & [6, 35] & [29, 80] \\ \hline
3	&	(100, 150] & 33.33\% & [8, 26] & [74, 132] \\ \hline
4	&	(150, 240]  & 31.86\% & [6, 27] & [125, 217] \\ \hline
	\end{tabular}%
	\label{durations}%
\end{table}% 

To generate the treatment duration for a patient, we follow the following procedure: We first draw a random number between 0 and 1 and use the associated probabilities in Table \ref{durations} in order to determine the patient's class. Next, pre-medication and infusion durations of this patient are determined by assuming a uniform distribution within the actual pre-medication and infusion intervals in Table \ref{durations} for her class. To justify the use of the uniform distribution for these times, we provide the chi-square goodness-of-fit test results and plots of fitted vs. actual durations in Appendix \ref{sec:actualvsfitted}. We generate a different set of durations for this patient in each scenario by taking into account her class. We then repeat this process for each patient in each instance. In generating these instances, we also aim to ensure that the total estimated treatment time results in an expected overtime around 60 minutes per nurse, in line with the requirements of the unit. For each instance, we set the overtime limit $L$ large enough to avoid infeasibility, but also small enough to avoid unfair nurse schedules. In our experiments, the overtime value for any nurse in any scenario does not exceed 147 minutes. Furthermore, none of our instances results in zero overtime, which is another indication that our instances are generated in line with the actual assignment of patients to half-shifts.
 
Our instances are labeled as $i\_j\_k$, where $i$, $j,$ and $k$ refer to the instance number, number of patients, and number of scenarios, respectively. Note that we provide the explicit composition of patient classes considered in each of the 10 instances in Appendix \ref{sec:instances}.

\subsection{Tests for varying CPLEX optimality gap value in the subproblems}
\label{testCPLEX}
A technique for solving scenario subproblems approximately based on optimality gap values that vary over iterations was introduced in Section \ref{gap}. The optimality gap values exhibit monotonically nonincreasing behavior over iterations due to the given rule. We perform experiments to see the effects of this technique on the computational times and quality of solutions by starting with two different initial optimality gap values of 10\% and 15\%. The results of the experiments are presented in Table \ref{tab:gapupdate}. In these experiments, the number of patients, nurses and chairs are set as 8, 2, and 4, respectively. 
\begin{table}
	\centering
	\caption{Effect of varying initial CPLEX optimality gap value on the objective value run time (in seconds), and number of iterations of the LPHA}
	\footnotesize
	%\renewcommand{\arraystretch}{1}
	%\centering\scalebox{0.9}{
		\begin{tabular}{|c|c|c|c|c|c|c|c|c|c|}
			\cline{2-10}
		    \multicolumn{1}{c|}{}& \multicolumn{3}{c|}{\textbf{Initial  gap = default}} & \multicolumn{3}{c|}{\textbf{Initial Gap = 0.10}} & \multicolumn{3}{c|}{\textbf{Initial Gap = 0.15}} \\
			\hline
			\textbf{Instance} & \textbf{Objective} & \textbf{Time} & \textbf{Iter} &  \textbf{Objective} & \textbf{Time} & \textbf{Iter} & \textbf{Objective} & \textbf{Time} & \textbf{Iter} \\
		\hline
    1\_8\_50 & 48.28 & 1491.97 & 78    & 48.46 & 2123.48 & 129   & 47.78 & 2357.88 & 128 \\
    \hline
    2\_8\_50 & 39.72 & 5444.85 & 107   & 46.82 & 3726.02 & 111   & 47.44      &   7092.87    & 121 \\
    \hline
    3\_8\_50 & 47.99 & 2467.34 & 96    & 60.61 & 2222.53 & 121   & 59.76 & 3117.67 & 113 \\
    \hline
    4\_8\_50 & 30.40 & 5137.10 & 127   & 34.95 & 2704.92 & 109   & 38.31 & 2964.76 & 172 \\
    \hline
    5\_8\_50 & 49.52 & 2037.69 & 103   & 54.03 & 3213.19 & 127   & 54.05 & 3623.33 & 117 \\
    \hline
    6\_8\_50 & 50.39 & 2266.95 & 120   & 49.05 & 3019.57 & 121   & 49.45 & 2705.13 & 118 \\
    \hline
    7\_8\_50 & 38.97 & 4881.50 & 103   & 45.65 & 5633.31 & 123   & 47.76 & 6575.99 & 141 \\
    \hline
    8\_8\_50 & 53.21 & 3140.59 & 105   & 53.07 & 2969.72 & 126   & 57.76 & 3111.89 & 112 \\
    \hline
    9\_8\_50 & 51.96 & 1449.22 & 108   & 57.64 & 2949.11 & 123   & 58.42 & 3519.49 & 111 \\
    \hline
    10\_8\_50 & 50.09 & 2262.28 & 94    & 48.46 & 4049.13 & 123   & 48.97 & 3730.82 & 116 \\
    \hline
    \textbf{Average} & \textbf{46.05} & \textbf{3057.95} & \textbf{104}   & \textbf{49.87} & \textbf{3173.54} & \textbf{121}   & \textbf{50.97} & \textbf{3879.98} & \textbf{125} \\
    \hline
	\end{tabular}%}
	\label{tab:gapupdate}%
\end{table}%

The results are compared with the version of the LPHA that uses default CPLEX optimality gap value while solving subproblems. There is no guarantee that the overall computation times decrease when the optimality gap is increased, because the algorithm may converge after larger amount of iterations. Indeed, the one that uses default CPLEX optimality gap value performs the best in terms of solution quality and run time on average. Due to insufficient basis for the use of rule, we do not vary the optimality gap in our further experiments. 

We also test various parameters of the LPHA to determine the best values to use in further experiments. The details of these experiments are provided Appendix \ref{sec:CR}.
 Consequently, we set the step size $\alpha$ as 2, penalty parameter $\rho$ as 0.0001, upper bound $\rho_1^u$ for the penalty parameter as 0.1, and the iteration limit as 100.

\subsection{Assessment of the LPHA Performance}
In this section, we initially evaluate the performance of the LPHA by conducting comparisons with optimal solutions. Next, the LPHA solutions are compared with the baseline schedule (i.e., actual schedule of the particular outpatient chemotherapy unit at Hacettepe Oncology Hospital). Finally, the LPHA is tested against commonly used scheduling heuristics from the literature.
 
%In this section, we fix the value of $\lambda$ at 0.3.
\subsubsection{Comparison with Optimal Solutions}
To validate the performance of the LPHA, we first compare our solutions to optimal solutions for instances that can be solved to optimality by CPLEX. Since CPLEX cannot find the optimal solutions in three hours of time limit even for problem instances having 8 patients, 5 scenarios, 2 nurses and 4 chairs, we decrease the number of patients to 7 for these experiments. We use five different $\lambda$ values on 10 different instances. 

On average, CPLEX spends approximately 100 seconds to find the optimal solution, whereas the LPHA spends approximately 40 seconds. The average optimality gap for the LPHA solutions is 5.64\%, which shows that the LPHA finds high quality solutions within significantly lower computational times than CPLEX.
%while it is less than 0.35\% in 25\% of the instances.
\begin{figure}
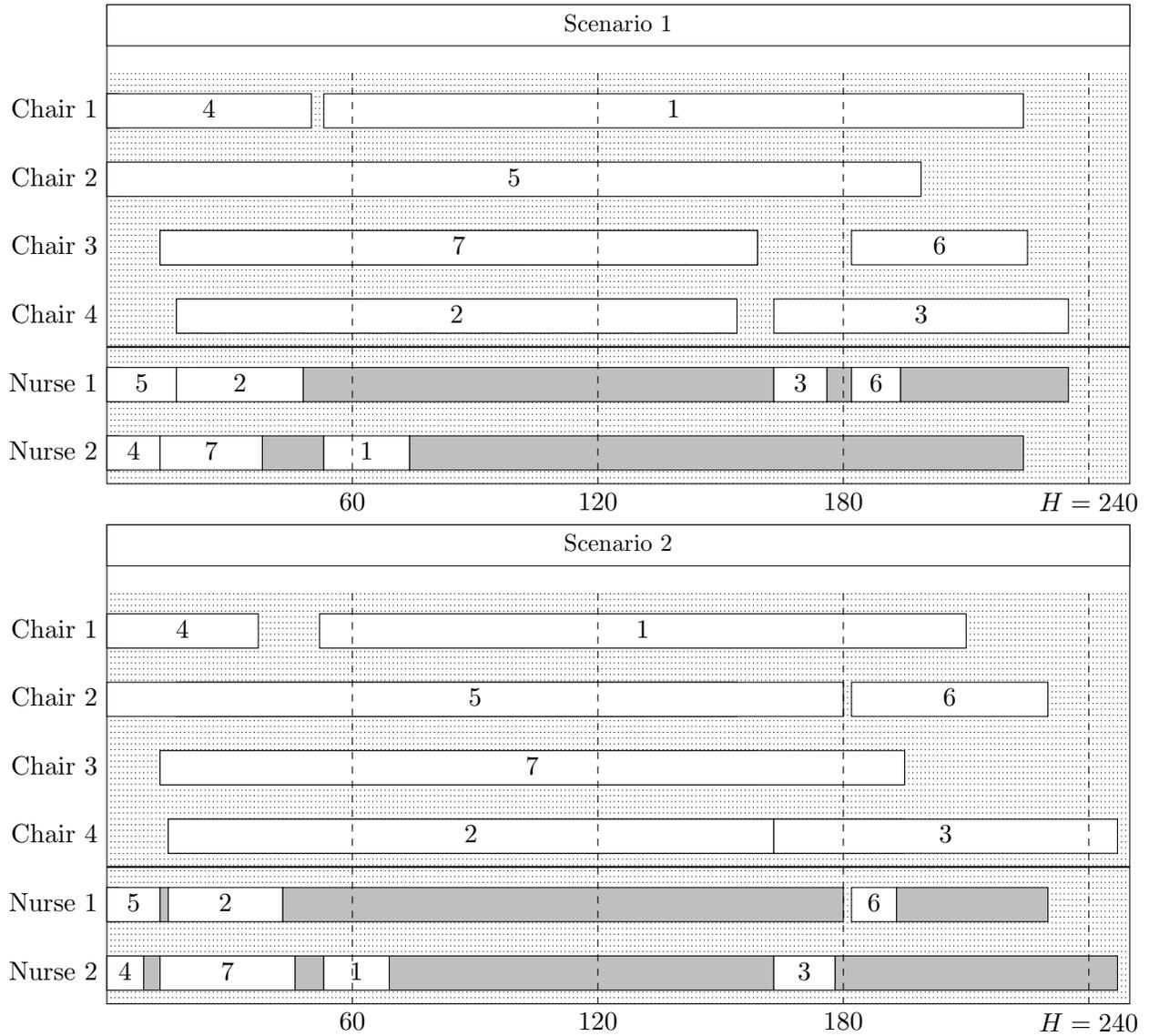

\begin{ganttchart}[vgrid,expand chart=\textwidth, bar height = 0.5]{1}{250}
\gantttitle{Scenario 1}{250}\\
\ganttbar{Chair 1}{1}{3}
\ganttbar[inline]{4}{1}{50}
\ganttbar[inline]{1}{54}{224}\\
\ganttbar{Chair 2}{1}{3} 
\ganttbar[inline]{5}{1}{199}\\
\ganttbar{Chair 3}{14}{159}
\ganttbar[inline]{7}{14}{159}
\ganttbar[inline]{6}{183}{225}\\
\ganttbar{Chair 4}{18}{154}
\ganttbar[inline]{2}{18}{154}
\ganttbar[inline]{3}{164}{235} \ganttnewline[thick]

\ganttbar{Nurse 1}{1}{3}
\ganttbar[inline]{5}{1}{17}
\ganttbar[inline]{2}{18}{48}
\ganttbar[inline, bar/.append style={fill=gray!50}]{}{49}{163}
\ganttbar[inline]{3}{164}{176}
\ganttbar[inline, bar/.append style={fill=gray!50}]{}{177}{182}
\ganttbar[inline]{6}{183}{194}
\ganttbar[inline, bar/.append style={fill=gray!50}]{}{195}{235}\\
\ganttbar{Nurse 2}{1}{3}
\ganttbar[inline]{4}{1}{13}
\ganttbar[inline]{7}{14}{38}
\ganttbar[inline, bar/.append style={fill=gray!50}]{}{39}{53}
\ganttbar[inline]{1}{54}{74}
\ganttbar[inline, bar/.append style={fill=gray!50}]{}{75}{224}

\ganttvrule[vrule/.append style={thin}]{60}{60}
\ganttvrule[vrule/.append style={thin}]{120}{120}
\ganttvrule[vrule/.append style={thin}]{180}{180}
\ganttvrule[vrule/.append style={thin}]{$H=240$}{240}
\end{ganttchart}

\begin{ganttchart}[vgrid,,expand chart=\textwidth, bar height = 0.5]{1}{250}
\gantttitle{Scenario 2}{250}\\
\ganttbar{Chair 1}{1}{3}
\ganttbar[inline]{4}{1}{37}
\ganttbar[inline]{1}{53}{210}\\
\ganttbar{Chair 2}{18}{154}
\ganttbar[inline]{5}{1}{180}
\ganttbar[inline]{6}{183}{230}\\
\ganttbar{Chair 3}{14}{159}
\ganttbar[inline]{7}{14}{195}\\
\ganttbar{Chair 4}{16}{195} 
\ganttbar[inline]{2}{16}{163}
\ganttbar[inline]{3}{164}{247}
\ganttnewline[thick]

\ganttbar{Nurse 1}{1}{3}
\ganttbar[inline]{5}{1}{13}
\ganttbar[inline, bar/.append style={fill=gray!50}]{}{14}{15}
\ganttbar[inline]{2}{16}{43}
\ganttbar[inline, bar/.append style={fill=gray!50}]{}{44}{180}
\ganttbar[inline]{6}{183}{193}
\ganttbar[inline, bar/.append style={fill=gray!50}]{}{194}{230}\\
\ganttbar{Nurse 2}{1}{3}
\ganttbar[inline]{4}{1}{9}
\ganttbar[inline, bar/.append style={fill=gray!50}]{}{10}{13}
\ganttbar[inline]{7}{14}{46}
\ganttbar[inline, bar/.append style={fill=gray!50}]{}{47}{53}
\ganttbar[inline]{1}{54}{69}
\ganttbar[inline, bar/.append style={fill=gray!50}]{}{70}{163}
\ganttbar[inline]{3}{164}{178}
\ganttbar[inline, bar/.append style={fill=gray!50}]{}{179}{247}
\ganttvrule[vrule/.append style={thin}]{60}{60}
\ganttvrule[vrule/.append style={thin}]{120}{120}
\ganttvrule[vrule/.append style={thin}]{180}{180}
\ganttvrule[vrule/.append style={thin}]{$H=240$}{240}
\end{ganttchart}
\caption{Optimal 4-hour chair and nurse schedules for two scenarios of Instance 1, where each number corresponds to a patient and the appointment times of patients are $a_1=52$, $a_2=16$, $a_3=163$, $a_4=a_5=0$, $a_6=182$, and $a_7=13$}
\label{fig:schedule}
\end{figure}

Using these instances, we are also able to observe the structure of optimal solutions. For this end, Figure \ref{fig:schedule} shows the optimal chair and nurse schedules for the first two scenarios of the first instance in this set. One common theme in both solutions is the tight scheduling of appointment times in the beginning of the shift until all chairs are occupied. The remaining patients are scheduled to arrive in sequence based on the expected departure times of the patients currently receiving treatment, with some buffer time to avoid waiting. Indeed, scenario 1 results in only 2 minutes of patient waiting (for patient 2), whereas there is no patient waiting in scenario 2. Scenario 1 also results in no overtime, as patients arriving towards the end of the shift need less time for treatment compared to scenario 2, where 7 minutes of overtime is needed for nurse 2 to complete the treatment of patient 3.

From Figure \ref{fig:schedule}, we also observe the importance of assigning chairs and nurses to patients upon arrival, rather than a priori. For example, fixing the chair assignment of patients 5, 6, and 7 over all possible scenarios would result in either waiting of patients or nurse overtime, unless more resources are added. These scenarios also reveal the importance of relaxing the slot-based scheduling structure, since hourly or 2-hourly slots would lead to more patient waiting time and/or more nurse overtime in both scenarios.

\subsubsection{Comparison with the Baseline Schedule}
To demonstrate the benefit of using LPHA to create schedules at the Outpatient Chemotherapy Unit at Hacettepe Oncology Hospital, we compare the LPHA solutions with the baseline schedule. Since the LPHA solves instances that represent half-day shifts, we generate two schedules using our algorithm for a representative day in our data set. The resulting morning and afternoon schedules are compared with the corresponding shifts of the baseline schedule for the same day on Table 3.   Note that we set $\lambda_{1}, \lambda_{2}$, $\lambda_{3}$ values as 0.3, 0.3, 0.4, respectively. The number of patients, nurses and chairs are fixed as 8, 2, and 4 in these experiments, respectively.

While selecting the day for schedule comparison, we consider the percentage of patient classes observed in our data set, which is given on Table \ref{durations}. We then identify a day in which the patient class composition is consistent with those average values. It should be noted here that the fixed-time slots used by the Outpatient Chemotherapy Unit start at 8:00 and 10:30 for the morning shift, and 13:00 and 15:30 for the afternoon shift. Since the chair set considered in each instance consists of four chairs and eight patients are scheduled for the set at each shift, four patients are asked to arrive at the beginning of each time slot in the baseline schedule. The LPHA schedule is also generated based on the same patient class composition. Then, using the same scenario set for treatment durations, the baseline schedule is compared with the LPHA schedule. To find the performance measures associated with the baseline schedule for each shift, the second-stage problem of the SMIP is solved.  

Table \ref{tab:baseline} compares the LPHA and baseline schedules with respect to average waiting time, nurse overtime and chair idle time. All three measures are improved by the LPHA at both shifts. The use of fixed-length slots at the Outpatient Chemotherapy Unit explains the reason behind the comprehensive improvement. The scope of improvement is particularly large with respect to average waiting time. The tendency to underestimate the majority of the appointment durations in the unit results in excessive wait amounts for patients. Note that the rate of reduction of nurse overtime is also notable, particularly for the morning shift. 

\begin{table}
	\centering
	\caption{Average patient waiting time, nurse overtime, and chair idle time of the LPHA and baseline schedules}
	\small{
	\begin{tabular}{|c|c|c|c|c|c|c|}
		\cline{2-7}    \multicolumn{1}{r|}{} & \multicolumn{3}{c|}{\textbf{Morning Shift}} & \multicolumn{3}{c|}{\textbf{Afternoon Shift}} \\
		\cline{2-7}    \multicolumn{1}{r|}{} & \textbf{Waiting time} & \textbf{Overtime} & \textbf{Idle time} & \textbf{Waiting time} & \textbf{Overtime} & \textbf{Idle time} \\
		\hline
		\textbf{LPHA} & 11.70  & 2.78  & 169.98 & 30.88 & 11.28 & 104.14 \\
		\hline
		\textbf{Baseline} & 86.66 & 18.84 & 185.58 & 146.94 & 15.70  & 110.82 \\
		\hline
	\end{tabular}}
	\label{tab:baseline}%
\end{table}%

\subsubsection{Comparison with Scheduling Heuristics}
\label{sec:JH}
To compare the performance of the LPHA with practice in the other outpatient chemotherapy units, we assess its performance against various combinations of sequencing and appointment time setting heuristics from the relevant literature. We sequence patients using four different sequencing heuristics: increasing mean of treatment time (SPT), decreasing mean of treatment time (LPT), increasing variance of treatment time (VAR), and increasing coefficient of variation of treatment time (CoV). We consider a \textit{job hedging} heuristic to determine the estimated lengths of patient appointments. 

Before discussing the implementation details of the heuristics, the structure of the model in terms of first and second-stage decisions is investigated. When values for patient precedence $b_{ij}$ and appointment time $a_{i}$ variables become known in the first-stage, it is then easy to settle on the second-stage decisions. In some appointment scheduling studies, such as Mancilla and Storer (2012)\nocite{Mancilla2012}, the second-stage subproblems are considered as network flow problems. Having multiple resources in our case makes the subproblems more challenging. However, the optimal nurse and chair assignment decisions at this stage can still be made easily using a simple rule that checks $b_{ij}$ values. The rule would favor the patient preceding others in the list while assigning the first available nurse and chair. Next, actual treatment start time for each patient can be determined accordingly, and their waiting times are revealed according to the discharge time of the previously scheduled patients. By checking discharge times of the patients treated by a nurse, overtime value for each nurse can be easily calculated. Since solving the second-stage problem is easy, finding a heuristic approach to solve the first-stage problem is a more critical issue. After setting the first-stage variable values, we call CPLEX rather then using the simple rule to solve the second-stage problem, as CPLEX finds the optimal solution within a very short amount of time.

To create schedules, we first choose one of the four sequencing heuristics discussed above. There are two different durations that may affect the sequencing rules, namely average pre-medication and infusion durations. We use average treatment time, which is equal to the summation of those two values. 

After the sequence is fixed, we utilize the job hedging heuristic to set patient appointment times (Yellig and Mackulak, 1997)\nocite{Yellig1997}. We apply the heuristic in accordance with the approach also used in Gul (2018)\nocite{gul2018stochastic}, Castaing et al. (2016)\nocite{Castaing2016}, and Gul et al. (2011)\nocite{gul2011bi}. We sort the durations for the relevant patient class in the data set in non-increasing order, and calculate the $k^{th}$ percentile of the set to use it for determining the estimated treatment duration of a patient. Since there are two different chemotherapy durations in the problem, the percentile values associated with them are separately determined, and the summation of those two values is considered while setting patient appointment times. The steps of the scheduling heuristic that combines job hedging with one of the sequencing heuristics, namely LPT, is demonstrated in Algorithm \ref{alg:JH}.
\begin{algorithm}
	\caption{Job Hedging Heuristics for LPT Rule}\label{alg:JH}
	\begin{algorithmic}[1]
		\State Calculate the average treatment duration for every patient using $avg_{i} = \frac{{\sum_{\omega \in \Omega} (s^{\omega}_{i} + t^{\omega}_{i})}}{|\Omega|} \quad \forall i \in I$.
		\State Sort the values of the $avg_{i}$ according to LPT rule, assign index numbers to patients and sequence the patients according to their index numbers.
		\State Assign $b_{ij}$ values according to index numbers found at Step 2.
		\State Set appointment times according to given percentile level of job hedging heuristic. If the index of the patient is less than or equal to the both of the numbers of chair and nurse, the appointment time of the associated patient is assigned as zero. The appointment times of the other patients are set sequentially according to their index numbers by checking the first estimated available time of nurse and chair simultaneously.
		\State Call CPLEX to obtain second-stage decisions and objective function value.
	\end{algorithmic}
\end{algorithm}

We vary the percentile values between 40\% and 65\% in increments of 5\% for all sequencing rules and fix the values of $\lambda_{1}, \lambda_{2},$ and $\lambda_{3}$ as 0.1, 0.8, and 0.1, respectively. Table 3 compares the performances of the heuristics with that of the LPHA by illustrating the average gap between the solutions of the LPHA and each combination of sequencing and job hedging heuristic for a specific percentile level.

\begin{table}
	\centering
	\caption{Average percent gaps of different sequencing heuristics from the LPHA solutions under varying job hedging levels}
	\begin{tabular}{|c|c|c|c|c|}
		\hline
		\textbf{Job hedging level} & \textbf{SPT gap} & \textbf{LPT gap} & \textbf{VAR gap} & \textbf{CoV gap} \\
		\hline
		40\%  & 80.1 & 23.5 & 71.3 & 57.2 \\
		\hline
		45\%  & 78.5 & 26.1 & 71.7 & 57.6 \\
		\hline
		50\%  & 83.2 & 32.6 & 70.5 & 60.5 \\
		\hline
		55\%  & 79.8 & 37.8 & 75.9 & 63.7 \\
		\hline
		60\%  & 89.7 & 48.3 & 80.2 & 69.5 \\
		\hline
		65\%  & 92.5 & 57.9 & 83.6 & 74.5 \\
		\hline
		\textbf{Average} & \textbf{84.0} & \textbf{37.7} & \textbf{75.5} & \textbf{63.8} \\
		\hline
	\end{tabular}
	\label{tab:heur}%
\end{table}%
The results show that the LPHA outperforms all combinations of heuristics significantly. On the average, LPT rule performs better than the other sequencing rules. This result makes sense since the patients with longer treatment durations are assigned to the earlier hours of the day also at the outpatient chemotherapy unit in the Hacettepe Oncology Hospital. Note that the LPHA improves the solutions of even the best-performing and commonly used rule by 37.7\%. Furthermore, the percentiles 45\%, 40\%, 50\%, and 40\% used for assigning appointment times are the best in terms of solution quality on average for sequencing rules SPT, LPT, VAR, and CoV, respectively. 

We also compare the LPHA with commonly used sequencing rules without using job hedging. The sequence of patients is fixed with one of the four aforementioned sequencing rules. Then, CPLEX is called to determine the appointment times and nurse/chair assignments. We call these modified versions of the heuristics as \textit{SPT-opt, LPT-opt, VAR-opt}, and \textit{CoV-opt}, respectively. Although we fix the sequence of the patients, assigning appointment times to the patients still remains a challenging optimization problem for CPLEX due to the uncertainty of the treatment durations. Therefore, we put a three-hour time limit on CPLEX and report the best solution for the given sequence. In Table \ref{tab:opt2}, we compare the results in terms of average gap between the solutions of the LPHA and each type of sequencing rules.

\begin{table}
	\centering
	\caption{Average percent gaps of different sequencing heuristics from the LPHA solutions}
	\begin{tabular}{|c|c|}
	\hline
	\textbf{Sequencing Rule} & \textbf{Average Gap} \\
	\hline
		SPT-opt     & 78.9 \\
		\hline
		LPT-opt   & 11.9 \\
		\hline
		VAR-opt   & 69.9 \\
		\hline
		CoV-opt    & 46.3 \\
		\hline
	\end{tabular}
	\label{tab:opt2}%
\end{table}%

As shown in Table \ref{tab:opt2}, the LPHA outperforms all heuristics, even when the appointment times and second-stage decisions are determined by a mathematical model, rather than simple job hedging rules. LPT-opt, which outperforms all other heuristics, performs 11.9\% worse on average than the LPHA in terms of solution quality.  

\begin{table}
\centering
\caption{Appointment times (in minutes) for the eight patients, expected waiting time (EWT), expected nurse overtime (EOT), expected chair idle time (EIT), and objective values under different heuristic approaches for the instance 1\_8\_50}
\label{tab:apptimes}
\small
\begin{tabular}{|c|c|c|c|c|c|c|c|}
\hline
\textbf{Patient} & \textbf{LPHA} & \textbf{LPT-opt} & \textbf{LPT+40\%} & \textbf{LPT+65\%} & \textbf{SPT-opt} & \textbf{SPT+40\%} & \textbf{SPT+65\%} \\
\hline
1	&	225	&	222	&	191	&	214	&	0	&	0	&	0	\\ \hline
2	&	17	&	11	&	14	&	20	&	52	&	34	&	40	\\ \hline
3	&	180	&	182	&	176	&	195	&	6	&	8	&	11	\\ \hline
4	&	148	&	148	&	170	&	195	&	8	&	8	&	11	\\ \hline
5	&	0	&	0	&	0	&	0	&	75	&	66	&	83	\\ \hline
6	&	0	&	17	&	16	&	22	&	39	&	33	&	40	\\ \hline
7	&	207	&	205	&	186	&	214	&	0	&	0	&	0	\\ \hline
8	&	0	&	0	&	0	&	0	&	121	&	134	&	147	\\ \hline
\textbf{EWT} & 51.48 & 50.44 & 76.36 & 53.26 & 38.76 & 43.16 & 18.14 \\ \hline
\textbf{EOT} & 54.20 & 57.88 & 65.38 & 87.64 & 96.48 & 87.90 & 93.56 \\ \hline
\textbf{EIT} & 69.66 & 74.56 & 74.82 & 139.90 & 82.26 & 78.10 & 92.00 \\ \hline
\textbf{Objective} & \textbf{55.47} & \textbf{58.80} & \textbf{67.42} & \textbf{89.43} & \textbf{81.29} & \textbf{82.45} & \textbf{85.86} \\
\hline
\end{tabular}
\end{table}

To further investigate the differences in the solutions between the LPHA and benchmark heuristics, Table \ref{tab:apptimes} provides the appointment times for the eight patients, expected waiting time, overtime, chair idle time, and objective value of the instance 1\_8\_50 under various heuristic approaches. Comparison of the LPHA appointment times to those under LPT-opt shows that the order of the patients in both cases is very similar and the differences between the appointment times are within 17 minutes. The LPHA sacrifices some waiting time to improve the nurse overtime, for which the objective weight is higher. The differences in the appointment times become more significant as job hedging heuristics are used and percentile values are increased, which explains the difference between the objective values. Under SPT, the order of patients is almost completely reversed from that under the LPHA. This results in shorter expected patient waiting times, but substantially increases the nurse overtime, leading to the inferiority of the SPT-based heuristics to the LPT-based counterparts and the LPHA.

Note that the relative performances of sequencing rules shown on Table \ref{tab:opt2} are consistent with the results of Table \ref{tab:heur}. Even though the average gaps improve when CPLEX is used to set appointment times instead of job hedging heuristics, the ranks of the sequencing rules do not change. Therefore, the LPT would be preferred by the managers creating schedules without using any sophisticated tool. However, the most important finding of this section is that the managers need a sophisticated approach rather than simple rules to determine patient sequences.

\subsection{Sensitivity Analysis on Model Parameters}
\label{sec:SA}
In this section, we present a sensitivity analysis on the model parameters to generate managerial insights. 
The value of $\lambda$ and number of nurses and chairs are varied for this purpose.
 
\subsubsection{Impact of the $\bm{\lambda}$ value}
The chemotherapy unit manager would typically consider the trade-off between patient waiting time, nurse overtime, and chair idle time while designing schedules. To illustrate the scope of this trade-off, we vary the $\lambda$ values, where ($\lambda_{1}$,$\lambda_{2}$,$\lambda_{3}$) are the objective function coefficients associated with patient waiting time, nurse overtime, and chair idle time, respectively. $\lambda$ values equal to zero refers to the extreme case where the related objective term is not important performance measure for the model. We test twelve set of different values of $\lambda$ on an instance set. We fix the  number of patients, nurses, and chairs to 8, 2, and 4, respectively in these experiments.

We set the $\lambda$ values systematically across experiments. While the weights for any two measures are kept the same, the weight of the third measure may be (i) one-tenth, (ii) one-half of the others, or (iii) twice, (iv) ten times as much as the others in an experiment. After fixing the weight of waiting time as 1, our weight-setting procedure results in the weight combinations shown in the first column of Table \ref{tab:lambdatable}. Note that we normalized the values of $\lambda$ in the experiments so that they could be consistent with the procedure in the other experiments.

\begin{table}
	\centering
	\caption{Average objective values (in minutes) and CPU times (in seconds) of the LPHA for different combinations of the objective weights}
		\begin{tabular}{|c|c|c|c|c|c|}
			\hline
			& \textbf{Objective} & \textbf{CPU} & \textbf{Waiting} & \textbf{Nurse} & \textbf{Chair} \\
			$\mathbf{(\lambda_1,\lambda_2,\lambda_3)}$ & \textbf{value} & \textbf{time} & \textbf{time} & \textbf{overtime} & \textbf{idle time} \\
			\hline
			(1, 10, 10) & 46.38 & 1919.54 & 142.45 & 45.23 & 37.93 \\
			\hline
			(1, 2, 2) & 55.73 & 2075.36 & 59.27 & 57.24 & 52.44 \\
			\hline
			(1, 0.5, 0.5) & 46.14 & 1742.22 & 20.99 & 67.03 & 75.54 \\
			\hline
			(1, 0.1, 0.1) & 20.86 & 1330.06 & 5.26  & 78.73 & 118.95 \\
			\hline
			(1, 0.1, 1) & 50.17 & 2553.53 & 27.32 & 78.02 & 70.24 \\
			\hline
			(1, 0.5, 1) & 54.33 & 2076.88 & 30.94 & 71.63 & 69.05 \\
			\hline
			(1, 2, 1) & 53.81 & 2181.50 & 47.13 & 54.74 & 58.62 \\	
			\hline
			(1, 10, 1) & 51.66 & 1924.10 & 90.05 & 47.41 & 55.83 \\
			\hline
			(1, 1, 0.1) & 40.83 & 2184.99 & 14.87 & 60.74 & 101.42 \\	
			\hline
			(1, 1, 0.5) & 49.29 & 1751.78 & 27.15 & 60.78 & 70.60 \\	
			\hline
			(1, 1, 2) & 54.17 & 1930.14 & 53.75 & 60.42 & 53.61 \\
			\hline
			(1, 1, 10) & 49.59 & 2310.34 & 111.78 & 60.02 & 42.33 \\
			\hline
	\end{tabular}
	\label{tab:lambdatable}%
\end{table}%

Table \ref{tab:lambdatable}  clearly illustrates the trade-off between patient waiting time, nurse overtime, and chair idle time with reference to $\lambda$ values. Average patient waiting time decreases significantly when only the relative value of $\lambda_{1}$ is increased. Average nurse overtime and chair idle time values also drop respectively as each of the $\lambda_{2}$ and $\lambda_{3}$ values is increased individually while keeping the remaining two constant. However, the impact of $\lambda_{2}$ and $\lambda_{3}$ on their respective performance measures is not as dramatic as the impact of  $\lambda_{1}$ on average waiting time.

When only the $\lambda_{1}$ value increases, both average nurse overtime and chair idle time increase. To minimize patient waiting time, the model assigns appointment times in wider time intervals, which leads to increased overtime and idle time values. This shows that there exists trade-off between patient waiting time and the sum of nurse overtime and chair idle time. When only the  $\lambda_{2}$ is increased, average waiting time increases, but average idle time decreases. This indicates that minimizing overtime also helps to minimize chair idle time, meaning that nurse overtime and chair idle time are not conflicting measures. When only the  $\lambda_{3}$ value is increased, average waiting time increases, while the average nurse overtime stays about the same. This implies that increasing  the $\lambda_{3}$ value alone may not always help improve average nurse overtime. Note that the cases where the $\lambda$ values are 0 and 1 are not realistic for the unit manager. The manager may choose a $\lambda$ value between those extreme points according to the goals of the units.

\subsubsection{Impact of the number of nurses and chairs}
The decision maker in the chemotherapy unit would also be interested in investigating the effect of using different numbers of chairs and nurses on the performance measures. For this purpose, we vary the number of nurses and chairs in our instance sets between 1 and 3, and 4 and 6, respectively. We ignore the case with 3 nurses and 4 chairs, since the corresponding nurse to chair ratio is not realistic. We report the average objective values, run times, and performance measures in the objective function based on 10 instances in Tables \ref{tab:nurse} and \ref{tab:nurse2}.

\begin{table}
	\centering
	\caption{Average objective value and CPU time (in seconds) under varying number of nurses and chairs}
	\begin{tabular}{|c|c|c|c|c|c|c|c|}
		\cline{3-8}    \multicolumn{1}{r}{} &       & \multicolumn{3}{c|}{\textbf{Average Objective Value}} & \multicolumn{3}{c|}{\textbf{Average CPU Time}} \\
		\cline{3-8}    \multicolumn{1}{r}{} &       & $\bm{|C|=4}$ & $\bm{|C|=5}$ & $\bm{|C|=6}$ & $\bm{|C|=4}$ & $\bm{|C|=5}$ & $\bm{|C|=6}$ \\
		\cline{2-8}
		\multicolumn{1}{r|}{} & $\bm{|N|=1}$ & 85.55 & 117.37 & 187.98 & 3052.16 & 2062.31 & 757.06 \\
		\cline{2-8}       \multicolumn{1}{r|}{} & $\bm{|N|=2}$ & 54.34 & 86.70 & 173.95 & \multicolumn{1}{c|}{1800.54} & 1469.18 & 776.65 \\
		\cline{2-8}        \multicolumn{1}{r|}{}  & $\bm{|N|=3}$ & Not solved & 82.82 & 172.42 & Not solved & 810.32 & 542.43 \\
		\cline{2-8}
	\end{tabular}%
	\label{tab:nurse}
\end{table}%

\begin{table}
	\centering
	\caption{Average waiting time, nurse overtime, and chair idle time (in minutes) under varying number of nurses and chairs}
	\footnotesize
	\begin{tabular}{|c|c|c|c|c|c|c|c|c|c|}
		\cline{2-10}  \multicolumn{1}{r}{}  & \multicolumn{3}{|c|}{\textbf{Waiting Time}} & \multicolumn{3}{c|}{\textbf{Nurse Overtime}} & \multicolumn{3}{c|}{\textbf{Chair Idle Time}} \\
		\cline{2-10}    \multicolumn{1}{r|}{} & $\bm{|C|=4}$ & $\bm{|C|=5}$ & $\bm{|C|=6}$ & $\bm{|C|=4}$ & $\bm{|C|=5}$ & $\bm{|C|=6}$ & $\bm{|C|=4}$ & $\bm{|C|=5}$ & $\bm{|C|=6}$ \\
		\hline
		$\bm{|N|=1}$ & 58.20 & 50.12 & 25.16 & 74.16 & 25.71 & 10.76 & 114.61 & 236.55 & 443.02 \\
		\hline
		$\bm{|N|=2}$ & 40.59 & 14.12 & 2.83  & 61.44 & 8.61  & 1.92  & 59.33 & 199.71 & 431.30 \\
		\hline
		$\bm{|N|=3}$ & Not solved & 8.52  & 0.72  & Not solved & 6.14  & 0.67  & Not solved & 196.06 & 430.01 \\
		\hline
	\end{tabular}
	\label{tab:nurse2}%
\end{table}%

As expected, when the number of nurses is increased, the average objective value decreases. All performance measures in the objective function (patient waiting time, nurse overtime, and chair idle time) are improved. The increase in the number of nurses provides to have more flexible schedules without extending clinic closure time. Therefore, both nurse overtime and chair idle time reduces. On the other hand, increase in the number of chairs results with having higher chair idle times. Therefore, the objective function value might increase depending on the net change in the decrease of patient waiting time and nurse overtime and increase in the chair idle time.
Another interesting issue is the change in the computation times. An increase in the number of chairs and nurses generally makes subproblems easier to solve resulting with a reduction in the solution times.

\subsection{Estimating the Value of Stochastic Solution (VSS)}
VSS is a cost of excluding uncertainty while giving the first-stage decisions. It measures the benefit of using stochastic programming solution over the mean value solution. We estimate the VSS by calculating the difference between the expected objective value associated with the mean value solution and the LPHA solution. In Table \ref{tab:VSS}, we report the relative VSS, which represents the percentage change in the objective values for various $\lambda$ values.

\begin{table}
	\centering
	\caption{Average objective values for the mean value solution (MV), LPHA solution, and relative VSS for varying objective weights}
	\begin{tabular}{|c|c|c|c|}
		\hline
		$\bm{(\lambda_{1},\lambda_{3},\lambda_{3})}$ & \textbf{MV} & \textbf{LPHA} & \textbf{\%VSS} \\
		\hline
		(0.3,0.3,0.4) & 58.56 & 54.34 & 7.7\% \\
		\hline
		(0.2,0.6,0.2) & 61.20 & 54.37 & 12.9\% \\
		\hline
		(0.8,0.1,0.1) & 39.67 & 25.12 & 58.8\% \\
		\hline
		(0.1,0.8,0.1) & 63.92 & 50.11 & 27.6\% \\
		\hline
		(0.1,0.1,0.8) & 66.24 & 50.46 & 32.7\% \\
		\hline
	\end{tabular}
	\label{tab:VSS}%
\end{table}%

As expected, VSS is always positive in our experimental runs. The LPHA improves the solutions of the mean value problem by 27.9\% on average. Therefore, considering uncertainty in chemotherapy appointment scheduling problem is valuable to minimize total weighted sum of patient waiting time, nurse overtime, and chair idle time. Note that the benefit of uncertainty modeling is particularly significant (with VSS\% of 58.8\%) for the cases where the highest priority is assigned to patient waiting time.

\section{Conclusion}
\label{chp:b6}
Over the recent years, due to the increase in the prevalence of cancer, the demand for chemotherapy units has been growing, and these units should have efficient planning and scheduling structures. In this paper, we focused on the Chemotherapy Appointment Scheduling Problem, where patients are sequenced and assigned appointment times by considering the availability of nurses and infusion chairs at the same time. The uncertainty in both the pre-medication and infusion durations are considered in the study. The aim is to design an appointment schedule in order to minimize the expected weighted sum of patient waiting time, nurse overtime, and chair idle time. A two-stage stochastic mixed integer programming formulation is used to formulate the problem.

The uncertain durations of chemotherapy are represented with scenarios in the stochastic programming formulation. Solving the problem even with a very few scenarios with CPLEX is computationally challenging. Therefore, a scenario decomposition based algorithm, namely the progressive hedging algorithm, is implemented to solve the problem. The PHA is enhanced through a number of modifications. A penalty parameter update method is proposed, where the parameters are set dynamically and controlled using changing limits.  A cycle detection method is used to guarantee convergence of the algorithm. A variable fixing procedure is incorporated to reduce overall computation times. Subproblem solution times are improved through symmetry-breaking constraints and bounds imposed on nurse overtime in the constraint set, and by linearization of the quadratic term in the objective function. The resulting algorithm after enhancements is called as linearized PHA (LPHA).

The proposed method is implemented based on real data from the Hacettepe Outpatient Chemotherapy Unit. To evaluate the performance of algorithm, the LPHA solutions are compared with the CPLEX solutions and the baseline schedule. Moreover, several scheduling heuristics are used to validate the performance of the algorithm. The results show that the LPHA outperforms commonly used heuristics in all cases. This finding implies that a sophisticated approach rather than simple rules is necessary to determine patient sequences in outpatient chemotherapy units. Furthermore, VSS is estimated to assess the value of considering uncertainty in our problem. It is found that the LPHA improves the solutions of mean value problem significantly. Our solution approach can be useful for chemotherapy unit managers to evaluate the effects of appointment schedules on nurse overtime,chair idle time, and patient waiting time. The unit managers can also observe the effects of changing the level of nurses and chairs in the chemotherapy unit. 

In a future study, the infusion schedules can be coordinated with lab test and oncologist evaluation schedules. A comprehensive model can be developed to determine treatment days and appointment times for patients simultaneously. Scenario subproblem solution routine in the LPHA can be parallelized to reduce computational time and solve larger size problem instances. Making the assumption of preserving patient arrival sequence throughout the whole chemotherapy unit visit allowed us to formulate the problem as a two-stage model. A multi-stage stochastic programming model can be formulated to study a variant of our problem by relaxing this assumption.

\bigskip
\bibliographystyle{plain}
\bibliography{thesis}
\bigskip
\newpage
\appendix

{\Large \textbf{Appendix}}

\section{Proof of Proposition 1} \label{sec:proof}

Let $\mathbf{B}=\{b_{ij},\, i,j \in I\}$ be the appointment precedence set from the first stage and let $\Phi_{F\!ACN}$ be the set of patient-chair and patient-nurse assignments (which we will simply call a \textit{schedule}) based on the first available chair and nurse assignment, given $\mathbf{B}$. Note that due to Constraints \eqref{bprec}, if $b_{ij}=1$, the appointment time of patient $i$ precedes that of patient $j$. We will prove the optimality of $\Phi_{F\!ACN}$ by showing that any optimal schedule $\Phi_{opt}$ can be converted to $\Phi_{F\!ACN}$ by a set of moves that result in no change in the objective function value.

In any step of the conversion from $\Phi_{opt}$ to $\Phi_{F\!ACN}$, let $\Phi_{cur}$ be the current schedule and let $i \in I$ be the earliest arriving patient in $\Phi_{cur}$ for whom the chair and/or nurse assignments are different than those in $\Phi_{F\!ACN}$. Let $C_1$ and $N_1$ denote the chair assignments for this patient under $\Phi_{cur}$, and $C_2$ and $N_2$ correspond to their counterparts under $\Phi_{F\!ACN}$. Three cases are possible.

\textbf{Case 1: $C_1 \neq C_2$, $N_1=N_2$}

This situation is illustrated in Figure \ref{fig:case1-after}(a). Let $k \in I$ be the first patient assigned to $C_2$ after the service of $i$ starts on $C_1$ and let $j$ be the last patient assigned to $C_2$ immediately preceding $k$. We denote the set of patient assignments before and after patient $i$ on $C_1$ by Subsequence $C_{11}$ and $C_{12}$, respectively. Furthermore, let Subsequence $C_{21}$ denote the set of assignments before $j$ and Subsequence $C_{22}$ those after $k$ on $C_2$. Note that here, patient $j$ must finish his/her service before that of $i$ starts on $C_1$, since otherwise $C_1$ would be available for $i$ before $C_2$, a contradiction to $C_2 \in \Phi_{F\!ACN}$ for patient $i$. Furthermore, due to Constraints \eqref{bprec}, $a_i \leq a_k$ must be preserved.

Next, consider swapping the assignments on $C_1$ for $i$ and Subsequence $C_{12}$ with those in $C_2$ for $k$ and Subsequence $C_{22}$, without changing the service start times and nurse assignments to patients (see Figure \ref{fig:case1-after}(b)). In such a case, total nurse overtime is constant (as nurse schedules remain the same), total patient waiting time is identical (as service start and end times are unchanged), and total chair idle time is as before (since the total makespan for $C_1$ and $C_2$ is equal to that under $\Phi_{cur}$ and total service time on these chairs is constant). Altogether, these imply that the overall objective function has stayed the same.

\begin{figure}
\begin{ganttchart}[inline,]{1}{30}
\ganttbar[inline=false]{$C_1$}{1}{12}
\ganttbar{Subsequence $C_{11}$}{1}{12}
\ganttbar{$i$}{13}{16}
\ganttbar{Subsequence $C_{12}$}{17}{25}
\\
\ganttbar{Subsequence $C_{21}$}{1}{9}
\ganttbar{$j$}{10}{11}
\ganttbar[inline=false]{$C_2$}{14}{19}
\ganttbar{$k$}{14}{19}
\ganttbar{Subsequence $C_{22}$}{20}{28}
\end{ganttchart}

\centering
(a)
\vspace{24pt}

\hspace{-0.8cm}
\begin{ganttchart}[inline,]{1}{30}
\ganttbar[inline=false]{$C_1$}{1}{12}
\ganttbar{Subsequence $C_{11}$}{1}{12}
\ganttbar{$k$}{14}{19}
\ganttbar{Subsequence $C_{22}$}{20}{28}
\\
\ganttbar{Subsequence $C_{21}$}{1}{9}
\ganttbar{$j$}{10}{11}
\ganttbar[inline=false]{$C_2$}{13}{16}
\ganttbar{$i$}{13}{16}
\ganttbar{Subsequence $C_{12}$}{17}{25}
\end{ganttchart}

\centering
(b)
\caption{The patient assignments for chairs $C_1$ and $C_2$ (a) under $\Phi_{cur}$ and (b) after swapping under Case 1 ($C_1 \neq C_2$, $N_1= N_2$)}
\label{fig:case1-after}
\end{figure}

\textbf{Case 2: $C_1 = C_2$, $N_1 \neq N_2$}

Since both chairs and nurses are identical among each other, we may use a symmetric approach to that in Case 1, exchanging the arguments for nurses and chairs. Starting with $m$ as the first patient assigned to nurse $N_2$ after the service of $i$ starts on $N_1$ (see Figure \ref{fig:case2-after}(a)), we perform a similar swap that results in a new set of assignments (see Figure \ref{fig:case2-after}(b)) with the chair assignments unchanged, where each patient starts his/her service at the same time as in $\Phi_{cur}$. This implies that the all three are still the same and objective function is identical to that under $\Phi_{cur}$.

\begin{figure}
\begin{ganttchart}[inline,]{1}{30}
\ganttbar[inline=false]{$N_1$}{1}{13}
\ganttbar{Subsequence $N_{11}$}{1}{13}
\ganttbar{$i$}{14}{17}
\ganttbar{Subsequence $N_{12}$}{18}{25}
\\
\ganttbar{Subsequence $N_{21}$}{1}{9}
\ganttbar[inline=false]{$N_2$}{15}{19}
\ganttbar{$l$}{10}{11}
\ganttbar{$m$}{15}{19}
\ganttbar{Subsequence $N_{22}$}{20}{28}
\end{ganttchart}

\centering
(a)
\vspace{24pt}

\hspace{-0.8cm}
\begin{ganttchart}[inline,]{1}{30}
\ganttbar[inline=false]{$N_1$}{1}{13}
\ganttbar{Subsequence $N_{11}$}{1}{13}
\ganttbar{$m$}{15}{19}
\ganttbar{Subsequence $N_{22}$}{20}{28}
\\
\ganttbar{Subsequence $N_{21}$}{1}{9}
\ganttbar[inline=false]{$N_2$}{14}{17}
\ganttbar{$l$}{10}{11}
\ganttbar{$i$}{14}{17}
\ganttbar{Subsequence $N_{12}$}{18}{25}
\end{ganttchart}

\centering
(b)

\caption{The patient assignments for nurses $N_1$ and $N_2$ (a) under $\Phi_{cur}$ and (b) after change under Case 2 ($C_1 = C_2$, $N_1 \neq N_2$)}
\label{fig:case2-after}
\end{figure}

\textbf{Case 3: $C_1 \neq C_2$, $N_1 \neq N_2$}

In this case, we perform a sequence of the Case 1 and Case 2 moves. By similar arguments, patient $i$ and the following subsequence assigned to $C_1$ can be swapped with patient $k$ and the following subsequence assigned to $C_2$, followed by swapping patient $i$ and the following subsequence assigned to $N_1$ with patient $m$ and the following subsequence assigned to $N_2$. By the above arguments, these two swaps are feasible and do not change the objective function value.

The conversion from $\Phi_{opt}$ to $\Phi_{F\!ACN}$ proceeds as follows: In each iteration, one checks the earliest arriving patient whose nurse and/or chair assignment is different from those in $\Phi_{F\!ACN}$. Depending on whether (i) chair, (ii) nurse, or (iii) both assignments are different, one applies the move in Case 1, 2, or 3, respectively, moving on to the next arriving patient. Since there are at most as many iterations as there are patients, the conversion can be made in no more than $|I|$ steps.
\QEDB

\section{SMIP-R}
\label{sec:DeterministicEquivalentModel}
\textbf{\underline{Revised Decision Variables:}} \\
$b^\omega_{ij}$  $=$  \(\begin{cases}
1, & \text{if patient $i \in I$ precedes patient $j \in I$ in daily appointment list in scenario $\omega \in \Omega$ } \\
0, & \text{otherwise} \end{cases}\)\\
$a^\omega_{i}$: appointment time of patient $i\in I$ in scenario $\omega \in \Omega$ \\
\footnotesize
	\begin{align}
	min \quad & \bigg(\lambda_{1} \displaystyle\sum_{i\in I}\displaystyle\sum_{\omega \in \Omega} p^\omega w^\omega_{i}+\lambda_{2} \displaystyle\sum_{n\in N}\displaystyle\sum_{\omega \in \Omega} p^\omega O^\omega_{n} +\lambda_{3} \displaystyle\sum_{c\in C}\displaystyle\sum_{\omega \in \Omega} p^\omega I^\omega_{c}  \bigg) & \\
	&\sum_{n\in N} x^\omega_{in} =1 \quad & \forall i \in I, \forall \omega \in \Omega \label{x2}\\
	&\displaystyle\sum_{c\in C} y^\omega_{ic} =1 \quad & \forall i \in I, \forall \omega \in \Omega \label{y2}\\
	&b^\omega_{ij} + b^\omega_{ji} = 1 \quad & \forall i, j \in I, j>i, \forall \omega \in \Omega \label{bprec3}\\
	&a^\omega_{i} + w^\omega_{i} + s^\omega_{i} + t^\omega_{i} = d^\omega_{i}  \quad & \forall i \in I, \forall \omega \in \Omega \label{timeeq2}\\
	&a^\omega_{j} + w^\omega_{j} \geq a^\omega_{i} + w^\omega_{i} + s^\omega_{i}-M(3-b_{ij}- x^\omega_{in}- x^\omega_{jn}) & \forall i, j \in I, j \neq i, \forall n \in N, \forall \omega \in \Omega \label{nurseprec2}\\ 
	&a_j + w^\omega_{j} \geq d^w_{i}-M(3-b_{ij}-y^\omega_{ic}-y^\omega_{jc}) \quad &\forall i, j \in I, j \neq i, \forall c \in C, \forall \omega \in \Omega \label{ypre2c} \\ 
	&a_j + w^\omega_{j} \geq a_i + w^\omega_{i} -M(1-b_{ij}) \quad & \forall i, j \in I, j \neq i, \forall \omega \in \Omega \label{bprec2}\\ 
	&a^\omega_j \geq a^\omega_i  -M(1-b^\omega_{ij}) \quad \quad & \forall i, j \in I, j \neq i, \forall \omega \in \Omega \\
	&O^\omega_{n} \geq d^\omega_{i} -H- M(1-x^\omega_{in}) \quad & \forall i \in I,\forall n \in N, \forall \omega \in \Omega \label{overtime2}\\
	&O^\omega_{n} \leq L \quad & \forall n \in N, \forall \omega \in \Omega \label{imp32}\\
	&T^\omega_{c} \geq H \quad & \forall c \in C, \forall \omega \in \Omega \\
	&T^\omega_{c} \geq d^\omega_{i}-M(1-y^{\omega}_{ic}) \quad & \forall i \in I, \forall c \in C, \forall \omega \in \Omega\\
	&I^\omega_{c} \geq T^\omega_{c}-\displaystyle\sum_{i\in I} (s^\omega_{i}+t^\omega_{i}) y^\omega_{ic} \quad & \forall c \in C, \forall \omega \in \Omega \\
	&a^\omega_{i} = a_{i} \quad & \forall i \in I, \forall \omega \in \Omega \label{nonant}\\
	& a^\omega_i :integer  \quad \quad & \forall i \in I, \forall \omega \in \Omega\\
	& b^\omega_{ij} \in \{0, 1\}  \quad \quad & \forall i, j \in I, \forall \omega \in \Omega \\
	&x^\omega_{in} \in \{0,1\}  \quad & \forall i \in I,\forall n \in N, \forall \omega \in \Omega \label{xsign2} \\
	&y^\omega_{ic} \in \{0,1\}  \quad & \forall i \in I,\forall c \in C, \forall \omega \in \Omega \label{ysign2} \\
	&d^\omega_{i},  w^\omega_{i} \geq 0  \quad & \forall i \in I, \forall \omega \in \Omega \label{dwsign2}\\
	&O^\omega_{n} \geq 0  \quad & \forall n \in N, \forall \omega \in \Omega \label{Osign2}\\
	&T^\omega_{c}, I^\omega_{c} \geq 0  \quad & \forall c \in C, \forall \omega \in \Omega
	\end{align} 

\normalsize

\section{Test Results and Distribution Plots of Actual and Fitted Data for Instance Generation} \label{sec:actualvsfitted}

In this section, we share the chi-square goodness-of-fit test results and comparisons of fitted vs. empirical distributions of pre-medication and infusion times for the four patient classes in Table \ref{durations}.

We first provide the goodness-of-fit results for the chi-square tests to justify the adherence of each pre-medication and infusion time distribution for each patient class to the given discrete uniform distribution in in Table \ref{durations}. For this purpose, we divide each time range into smaller categories, ensuring that the expected number of observations in each category based on the given fitted distribution is at least 5. Using a confidence level of 95\%, we test the null hypothesis that the measurements from our observation follows a discrete uniform distribution.

\begin{figure}
\centering
		\includegraphics[width=0.75\linewidth]{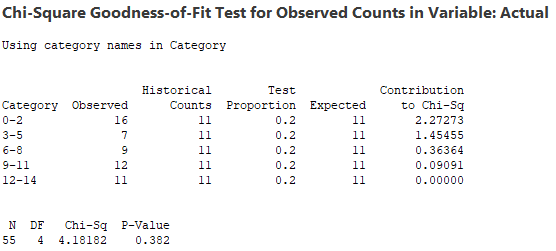}
		\caption{Chi-square goodness-of-fit results for the observed values of pre-medication times of patient class 1 to a discrete uniform distribution with parameters 0 and 14}
		\label{fig:gof1}
\end{figure}

\begin{figure}
\centering
		\includegraphics[width=0.75\linewidth]{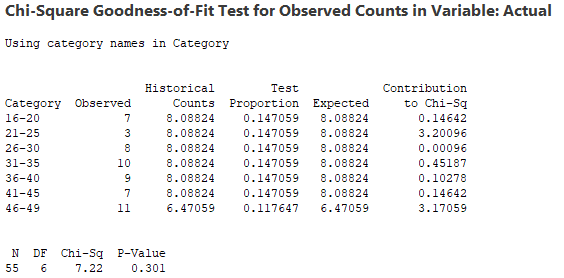}
		\caption{Chi-square goodness-of-fit results for the observed values of infusion times of patient class 1 to a discrete uniform distribution with parameters 16 and 49}
		\label{fig:gof2}
\end{figure}

\begin{figure}
\centering
		\includegraphics[width=0.75\linewidth]{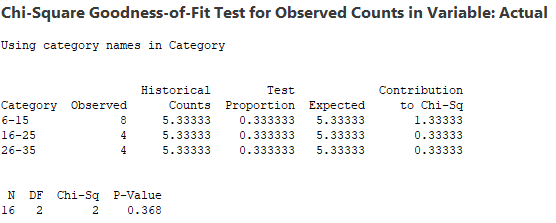}
		\caption{Chi-square goodness-of-fit results for the observed values of pre-medication times of patient class 2 to a discrete uniform distribution with parameters 6 and 35}
		\label{fig:gof3}
\end{figure}

\begin{figure}
\centering
		\includegraphics[width=0.75\linewidth]{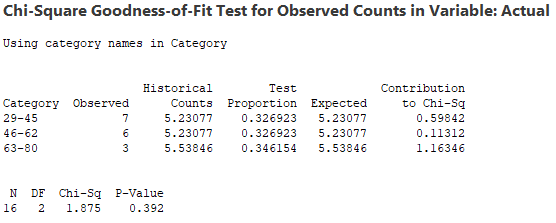}
		\caption{Chi-square goodness-of-fit results for the observed values of infusion times of patient class 2 to a discrete uniform distribution with parameters 29 and 80}
		\label{fig:gof4}
\end{figure}

\begin{figure}
\centering
		\includegraphics[width=0.75\linewidth]{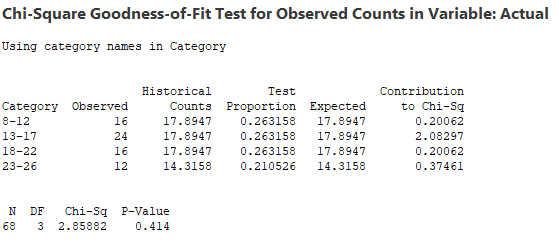}
		\caption{Chi-square goodness-of-fit results for the observed values of pre-medication times of patient class 3 to a discrete uniform distribution with parameters 8 and 26}
		\label{fig:gof5}
\end{figure}

\begin{figure}
\centering
		\includegraphics[width=0.75\linewidth]{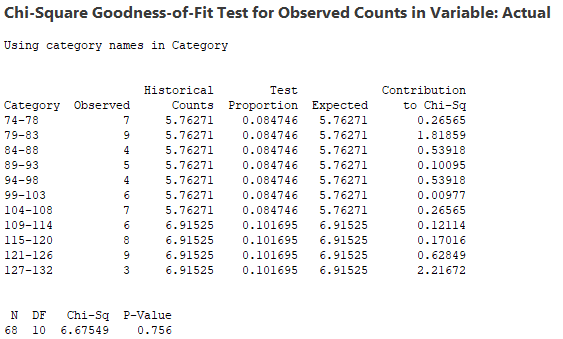}
		\caption{Chi-square goodness-of-fit results for the observed values of infusion times of patient class 3 to a discrete uniform distribution with parameters 74 and 132}
		\label{fig:gof6}
\end{figure}

\begin{figure}
\centering
		\includegraphics[width=0.75\linewidth]{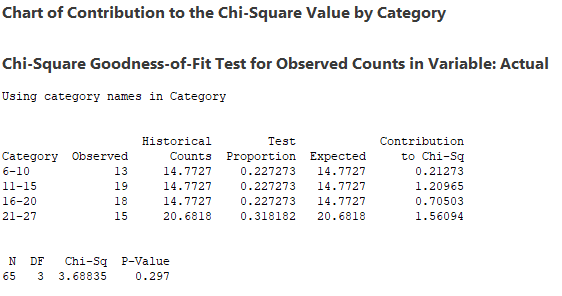}
		\caption{Chi-square goodness-of-fit results for the observed values of pre-medication times of patient class 4 to a discrete uniform distribution with parameters 6 and 27}
		\label{fig:gof7}
\end{figure}

\begin{figure}
\centering
		\includegraphics[width=0.75\linewidth]{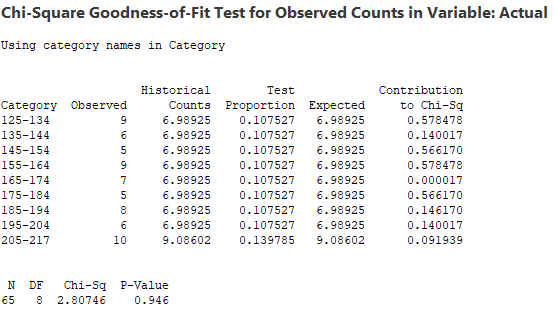}
		\caption{Chi-square goodness-of-fit results for the observed values of infusion times of patient class 4 to a discrete uniform distribution with parameters 125 and 217}
		\label{fig:gof8}
\end{figure}

Figures \ref{fig:gof1} through \ref{fig:gof8} present the chi-square goodness-of-fit results, in each of which the p-value is large enough so that the null hypotheses cannot be rejected, implying that the observed times fit the proposed discrete uniform distributions.

\begin{figure}
\centering
		\includegraphics[width=0.75\linewidth]{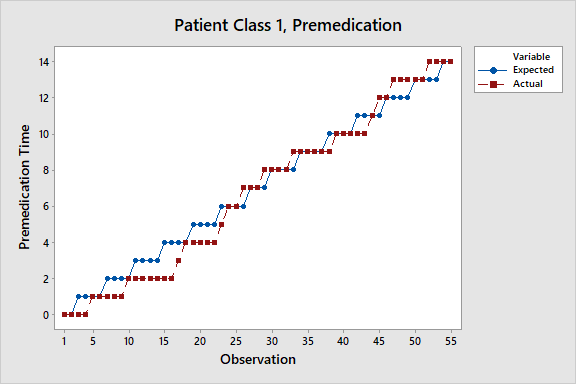}
		\caption{Fitted vs. actual pre-medication times of patient class 1}
		\label{fig:plot1}
\end{figure}

\begin{figure}
\centering
		\includegraphics[width=0.75\linewidth]{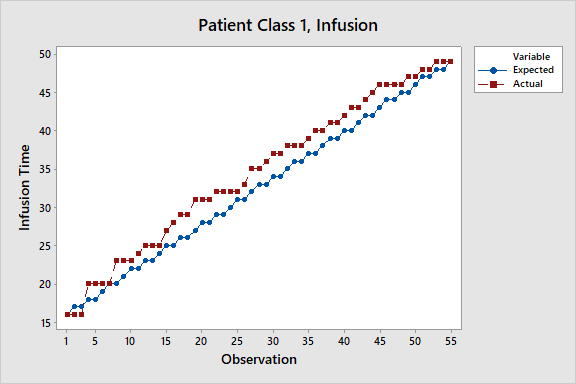}
		\caption{Fitted vs. actual infusion times of patient class 1}
		\label{fig:plot2}
\end{figure}

\begin{figure}
\centering
		\includegraphics[width=0.75\linewidth]{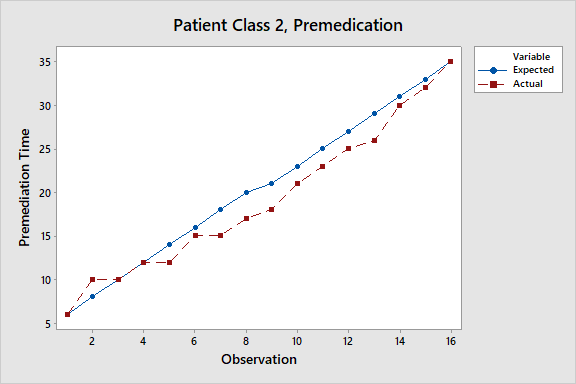}
		\caption{Fitted vs. actual pre-medication times of patient class 2}
		\label{fig:plot3}
\end{figure}

\begin{figure}
\centering
		\includegraphics[width=0.75\linewidth]{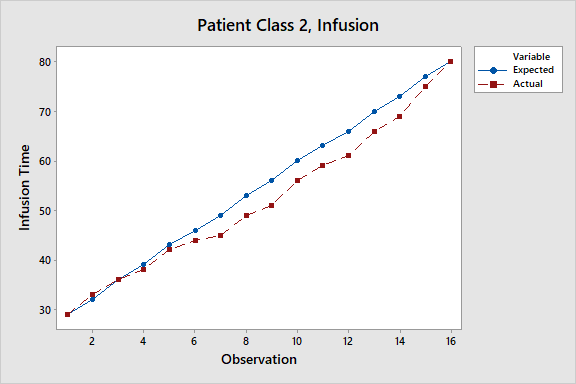}
		\caption{Fitted vs. actual infusion times of patient class 2}
		\label{fig:plot4}
\end{figure}

\begin{figure}
\centering
		\includegraphics[width=0.75\linewidth]{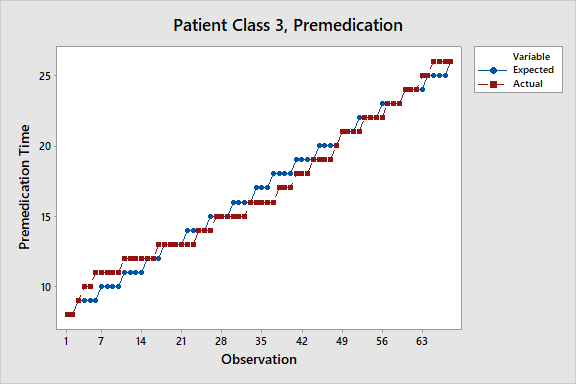}
		\caption{Fitted vs. actual pre-medication times of patient class 3}
		\label{fig:plot5}
\end{figure}

\begin{figure}
\centering
		\includegraphics[width=0.75\linewidth]{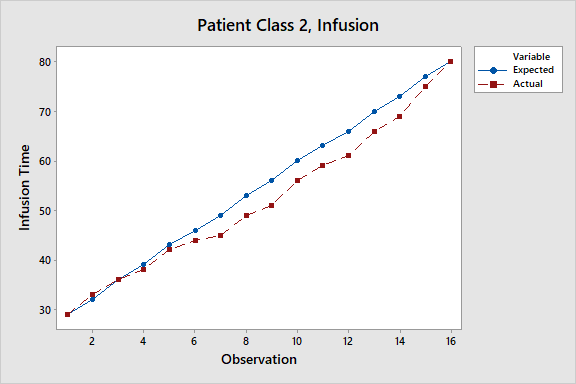}
		\caption{Fitted vs. actual infusion times of patient class 3}
		\label{fig:plot6}
\end{figure}

\begin{figure}
\centering
		\includegraphics[width=0.75\linewidth]{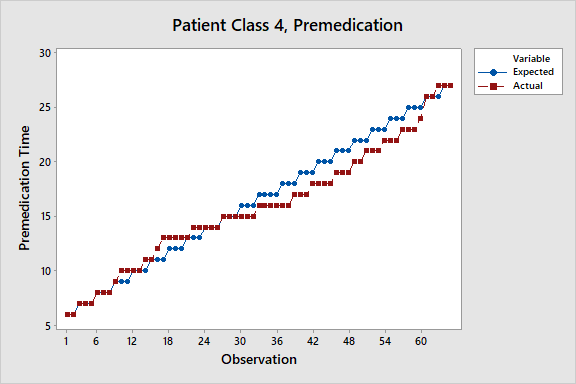}
		\caption{Fitted vs. actual pre-medication times of patient class 4}
		\label{fig:plot7}
\end{figure}

\begin{figure}
\centering
		\includegraphics[width=0.75\linewidth]{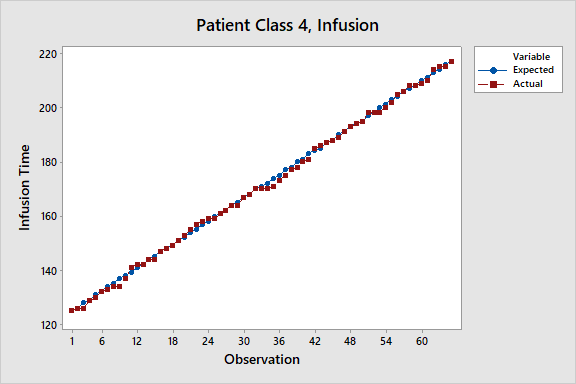}
		\caption{Fitted vs. actual infusion times of patient class 4}
		\label{fig:plot8}
\end{figure}

To provide further justification for fitting a discrete uniform distribution for the pre-medication and infusion times for each patient class, we provide plots of fitted and actual values in Figures \ref{fig:plot1} through \ref{fig:plot8}. As can be observed from these figures, the fitted values do not diverge significantly from the actual times in any case, which further allows us to use these distributions to generate the pre-medication and infusion times.

\newpage
\section{Patient Composition Details} \label{sec:instances}

Table \ref{tab:compositon} shows the percentage of each class of patients in the 10 instances used within the experiments.

\begin{table}[htb]
	\centering
	\caption{Percentage of patients in each class in each instance}
	\label{tab:compositon}
	\begin{tabular}{|c|c|c|c|c|}
		\hline
		\textbf{Instance} & \textbf{Class 1} & \textbf{Class 2} & \textbf{Class 3} & \textbf{Class 4} \\
		\hline
		1\_8\_50  & 25.0 & 12.5 & 12.5 & 50.0 \\
		\hline
		2\_8\_50  & 12.5 & 37.5 & 0.0 & 50.0 \\
		\hline
		3\_8\_50  & 12.5 & 25.0 & 12.5 & 50.0 \\
		\hline
		4\_8\_50  & 25.0 & 12.5 & 25.0 & 37.5 \\
		\hline
		5\_8\_50  & 12.5 & 25.0 & 12.5 & 50.0 \\
		\hline
		6\_8\_50  & 12.5 & 12.5 & 37.5 & 37.5 \\
		\hline
		7\_8\_50  & 12.5 & 25.0 & 25.0 & 37.5 \\
		\hline
		8\_8\_50  & 12.5 & 12.5 & 37.5 & 37.5 \\
		\hline
		9\_8\_50  & 12.5 & 25.0 & 12.5 & 50.0 \\
		\hline
		10\_8\_50  & 25.0 & 12.5 & 12.5 & 50.0 \\
		\hline
	\end{tabular}
\end{table}

\section{Sensitivity analysis on the LPHA parameters}
\label{sec:CR}
In this part, we conduct experiments to evaluate the sensitivity of the LPHA to the algorithm parameters. As there are many parameters associated with the LPHA, we perform preliminary experiments to fix some of these before elaborating on others. As a result, initial values of Lagrangian multipliers ($\mu^{\omega (v)}_{i}$) are fixed to zero, since there is no significant effect of changing this parameter. We also fix the number of patients, nurses, and chairs as 8, 2, and 5, respectively.
%in the experiments presented at this section to see the effect of the LPHA parameters clearly. In the experiments presented at this section, we assumed two nurses, four chairs, and eight patients are present in the unit. The value of $\lambda$ is also fixed to 0.3.

The effects of different values set for $\alpha$, $\rho$, $\rho^u_1$, and  $iterlimit$ are tested based on a one-way sensitivity analysis approach. First, we investigate the effect of $\alpha$ on the performance of the algorithm. 

\subsection{Effect of $\boldsymbol{\alpha}$ on the performance of the LPHA}
The parameter $\alpha$ is the step size that determines the rate of change in the penalty parameter $\rho$. The parameter $\rho$ may be updated by a factor of $\alpha$, 1, or $\dfrac{1}{\alpha}$ according to the changes in the $\Delta_{p}$ and $\Delta_{d}$, defined in Section 4. When $\alpha$ takes higher values, this means the penalty parameter varies significantly during the consecutive iterations. On the other hand, small values of $\alpha$ smoothly change the value of penalty parameter. When the deviation in the penalty parameter is smaller, we expect convergence in longer computational times to a higher quality solution. 

\begin{table}
	\centering
\caption{Effect of $\alpha$ on the objective value (Obj), run time (in seconds), and number of iterations (Iter) of the LPHA}
\begin{tabular}{|c|c|c|c|c|c|c|c|c|c|}
		\cline{2-10}    \multicolumn{1}{r|}{} & \multicolumn{3}{c|}{$\bm{\alpha=2}$} & \multicolumn{3}{c|}{$\bm{\alpha=4}$} & \multicolumn{3}{c|}{$\bm{\alpha=6}$} \\
		\hline
		\textbf{Instance} & \textbf{Obj} & \textbf{Time} & \textbf{Iter} & \textbf{Obj} & \textbf{Time} & \textbf{Iter} & \textbf{Obj} & \textbf{Time} & \textbf{Iter} \\
		\hline
		1\_8\_50 & 86.41 & 894.67 & 131   & 86.03 & 1025.06 & 155   & 86.27 & 1423.75 & 129 \\
		\hline
		2\_8\_50 & 97.01 & 1213.24 & 115   & 97.88 & 661.33 & 115   & 98.62 & 1025.54 & 116 \\
		\hline
		3\_8\_50 & 82.99 & 1978.06 & 142   & 83.31 & 1171.81 & 135   & 84.79 & 1662.01 & 118 \\
		\hline
		4\_8\_50 & 101.35 & 961.67 & 98    & 101.41 & 5070.40 & 122   & 102.98 & 1527.77 & 197 \\
		\hline
		5\_8\_50 & 83.25 & 1366.29 & 114   & 82.90 & 1748.01 & 199   & 81.14 & 1895.23 & 181 \\
		\hline
		6\_8\_50 & 86.49 & 2457.74 & 236   & 81.49 & 1004.17 & 120   & 82.32 & 1713.64 & 159 \\
		\hline
		7\_8\_50 & 93.39 & 1363.62 & 114   & 93.20 & 733.42 & 113   & 93.17 & 1416.79 & 117 \\
		\hline
		8\_8\_50 & 79.90 & 1786.17 & 166   & 80.40 & 803.35 & 118   & 84.85 & 1702.69 & 121 \\
		\hline
		9\_8\_50 & 79.23 & 1836.90 & 131   & 80.39 & 1144.87 & 127   & 79.56 & 2076.86 & 123 \\
		\hline
		10\_8\_50 & 80.78 & 1851.17 & 147   & 80.01 & 1329.35 & 164   & 80.99 & 1732.41 & 118 \\
		\hline
		\textbf{Average} & 87.08 & 1570.95 & 139.40 & 86.70 & 1469.18 & 136.80 & 87.47 & 1617.67 & 137.90 \\
		\hline
	\end{tabular}
	\label{tab:alphachange}
\end{table}

We set three different values of $\alpha$ (2, 4, and 6) and report their impact into computational times and solution quality in Table \ref{tab:alphachange}. In doing so, the aim is to fix $\alpha$ at a reasonable value where there is a balance between solution time and quality. Based on the results in Table \ref{tab:alphachange}, we observe that $\alpha=4$ results in shorter computational times than the other two cases (with CPU time savings of 6.9\% and 10.1\% over than of $\alpha=2$ and $\alpha=6$, respectively). Hence, for the remainder of the experiments, we fix $\alpha$ at this value.

\subsection{Effect of initial value of $\bm{\rho}$ on the performance of the LPHA}
There are various factors that affect the behaviour of the penalty parameter value $\rho$. The first one is its initial value. According to the literature, the lower values of $\rho$ yield better quality solutions in longer computational times. To test the effect of the initial value of $\rho$, three different values are used (0.0001, 0.005, and 0.1). The experimental results are provided in Table \ref{tab:rho}.

\begin{table}
	\centering
	\caption{Effect of the initial value of penalty parameter ($\rho$) on the objective value (Obj), run time (in seconds), and number of iterations (Iter) of the LPHA}
\begin{tabular}{|c|c|c|c|c|c|c|c|c|c|}
	\cline{2-10}    \multicolumn{1}{r|}{} & \multicolumn{3}{c|}{$\bm{\rho=0.0001}$} & \multicolumn{3}{c|}{$\bm{\rho=0.0050}$} & \multicolumn{3}{c|}{$\bm{\rho=0.1000}$} \\
		\hline
		\textbf{Instance} & \textbf{Obj} & \textbf{Time} & \textbf{Iter} & \textbf{Obj} & \textbf{Time} & \textbf{Iter} & \textbf{Obj} & \textbf{Time} & \textbf{Iter} \\
	\hline
	1\_8\_50 & 86.03 & 1025.06 & 155   & 86.99 & 1372.51 & 120 & 86.48 & 1081.56 & 151 \\
	\hline
	2\_8\_50 & 97.88 & 661.33 & 115   & 98.09 & 1018.9 & 114 & 98.89 & 771.58 & 116 \\
	\hline
	3\_8\_50 & 83.31 & 1171.81 & 135   & 88.77 & 1042.49 & 125 & 82.93 & 1436.79 & 159 \\
	\hline
	4\_8\_50 & 101.41 & 5070.40 & 122   & 105.30 & 711.86 & 108 & 102.45 & 823.04 & 120 \\
	\hline
	5\_8\_50 & 82.90 & 1748.01 & 199   & 81.76 & 1313.61 & 125 & 82.83 & 1588.85 & 193 \\
	\hline
	6\_8\_50 & 81.49 & 1004.17 & 120   & 81.13 & 1334.77 & 123 & 81.33 & 1001.15 & 154 \\
	\hline
	7\_8\_50 & 93.20 & 733.42 & 113   & 96.54 & 1057.86 & 126 & 93.42 & 897.90 & 121 \\
	\hline
	8\_8\_50 & 80.40 & 803.35 & 118   & 81.97 & 960.78 & 127 & 86.95 & 1131.55 & 143 \\
	\hline
	9\_8\_50 & 80.39 & 1144.87 & 127   & 81.62 & 1187.91 & 113 & 85.98 & 1319.30 & 106 \\
	\hline
	10\_8\_50 & 80.01 & 1329.35 & 164   & 82.77 & 1376.28 & 152 & 81.54 & 1304.55 & 154 \\
	\hline
	\textbf{Average} & \textbf{86.70} & \textbf{1469.18} & \textbf{136.80} & \textbf{88.50} & \textbf{1137.70} & \textbf{123.30} & \textbf{88.28} & \textbf{1135.63} & \textbf{141.70} \\
	\hline
\end{tabular}
\label{tab:rho}%
\end{table}%

The computational times decrease significantly in the cases where $\rho$ is equal to 0.005 and 0.1 compared to $\rho=0.0001$ case. Since the increase in the initial value of $\rho$ may create an undesirable increase in the objective values for some instances, we use an initial $\rho$ value of 0.0001 in the subsequent experiments. 
%The lower values of the penalty parameters provide consistently good objective function values. 
\subsection{Effect of $\bm{\rho^u_1}$ on the performance of the LPHA}
As stated in Section 4, the value of the penalty parameter $\rho$ is limited using an upper bound $\rho^u_1$. This prevents faster convergence of the algorithm to premature solutions due to high penalty parameter values. There is a trade-off between solution quality and computational times to be considered while selecting the value of $\rho^u_1$. For this purpose, we used three different values as 0.1, 0.5, and 0.7 to test the performance of the algorithm. The experimental results are provided in Table \ref{tab:pulimit}.

\begin{table}
	\centering
\caption{Effect of the first upper limit of penalty parameter ($\rho_1^u$) on the objective value (Obj), run time (in seconds), and number of iterations (Iter) of the LPHA}
	\begin{tabular}{|c|c|c|c|c|c|c|c|c|c|}
		\cline{2-10}    \multicolumn{1}{r|}{} & \multicolumn{3}{c|}{$\bm{\rho^u_1}=0.1$} & \multicolumn{3}{c|}{$\bm{\rho^u_1}=0.5$} & \multicolumn{3}{c|}{$\bm{\rho^u_1}=0.7$} \\
		\hline
		\textbf{Instance} & \textbf{Obj} & \textbf{Time} & \textbf{Iter} & \textbf{Obj} & \textbf{Time} & \textbf{Iter} & \textbf{Obj} & \textbf{Time} & \textbf{Iter} \\
		\hline
		1\_8\_50 & 86.03 & 1025.06 & 155   & 99.81 & 955.68 & 147   & 102.39 & 694.85 & 109 \\
		\hline
		2\_8\_50 & 97.88 & 661.33 & 115   & 98.61 & 925.16 & 169   & 99.62 & 682.63 & 131 \\
		\hline
		3\_8\_50 & 83.31 & 1171.81 & 135   & 81.11 & 964.57 & 131   & 81.58 & 943.45 & 147 \\
		\hline
		4\_8\_50 & 101.41 & 5070.40 & 122   & 104.66 & 1009.34 & 160   & 103.34 & 881.37 & 186 \\
		\hline
		5\_8\_50 & 82.90 & 1748.01 & 199   & 83.61 & 999.16 & 132   & 84.22 & 837.48 & 138 \\
		\hline
		6\_8\_50 & 81.49 & 1004.17 & 120   & 83.14 & 953.42 & 137   & 85.00 & 779.78 & 136 \\
		\hline
		7\_8\_50 & 93.20 & 733.42 & 113   & 93.09 & 1073.81 & 144   & 95.27 & 1388.00 & 150 \\
		\hline
		8\_8\_50 & 80.40 & 803.35 & 118   & 81.68 & 1070.32 & 126   & 84.52 & 1313.22 & 156 \\
		\hline
		9\_8\_50 & 80.39 & 1144.87 & 127   & 81.14 & 1028.52 & 125   & 81.59 & 752.32 & 115 \\
		\hline
		10\_8\_50 & 80.01 & 1329.35 & 164   & 80.62 & 1133.97 & 120   & 81.51 & 1001.34 & 122 \\
		\hline
		\textbf{Average} & \textbf{86.70} & \textbf{1469.18} & \textbf{136.80} & \textbf{88.75} & \textbf{1011.40} & \textbf{139.10} & \textbf{89.90} & \textbf{927.44} & \textbf{139.00} \\
		\hline
	\end{tabular}
	\label{tab:pulimit}%
\end{table}%

According to the results provided in Table \ref{tab:pulimit}, the best solutions with respect to the solution quality are obtained with the smallest tested value of $\rho^u_1$. The average number of iterations needed for convergence does not change significantly by varying the values of $\rho^u_1$. Since larger penalty parameter values make scenario subproblems easier to solve, the time spent for completing an iteration is expected to decrease. This explains the inverse relation between the value of $\rho^u_1$ and run time in the table.
When $\rho^u_1$ is raised to 0.5 and 0.7, the objective value deteriorates by 2.4\% and 3.7\% with respect to the case where $\rho^u_1=0.1$, respectively. Since the use of $\rho^u_1=0.1$ only needs, on average, less than ten extra minutes to improve objective value, we set 0.1 for $\rho^u_1$ for further experiments.

\subsection{Effect of %iteration limit for changing $\boldsymbol{\rho^u_1}$ to $\boldsymbol{\rho^u_2}$ 
\textit{iterlimit} on the performance of the LPHA}
In the LPHA, we use $\rho^u_1$ as an upper bound for the value of the penalty parameter up to a predetermined iteration $iterlimit$ in Algorithm 1. After $iterlimit$ many iterations, the upper bound on penalty limit is increased to $\rho^u_2$. We set reasonable values for $iterlimit$, because very large values of $iterlimit$ may lead to limited values of $\rho$ even when the iteration number reaches a considerable amount. This in turn may result in excessive computation times until convergence. On the other hand, having a limited value of $\rho$ over a large amount of iterations would improve solution quality. We try three different alternatives for determining the value of $iterlimit$ as 50, 70 and 100, because the computation times increase significantly once the $iterlimit$ exceeds 100, according to our preliminary experiments.

\begin{table} [htb]
	\centering
	\caption{Effect of iteration limit ($iterlimit$) on the objective value (Obj), run time (in seconds), and number of iterations (Iter) of the LPHA}
	\begin{tabular}{|c|c|c|c|c|c|c|c|c|c|}
		\cline{2-10}    \multicolumn{1}{r|}{} & \multicolumn{3}{c|}{$\bm{iterlimit}=50$} & \multicolumn{3}{c|}{$\bm{iterlimit}=70$} & \multicolumn{3}{c|}{$\bm{iterlimit}=100$} \\
		\hline
		\textbf{Instance} & \textbf{Obj} & \textbf{Time} & \textbf{Iter} & \textbf{Obj} & \textbf{Time} & \textbf{Iter} & \textbf{Obj} & \textbf{Time} & \textbf{Iter} \\
		\hline
		1\_8\_50 & 87.06 & \multicolumn{1}{c|}{604.769} & \multicolumn{1}{c|}{100} & 86.84 & \multicolumn{1}{c|}{733.507} & \multicolumn{1}{c|}{117} & 86.03 & 1025.06 & 155 \\
		\hline
		2\_8\_50 & 97.95 & \multicolumn{1}{c|}{401.448} & \multicolumn{1}{c|}{79} & 98.12 & \multicolumn{1}{c|}{445.828} & \multicolumn{1}{c|}{82} & 97.88 & 661.33 & 115 \\
		\hline
		3\_8\_50 & 81.56 & \multicolumn{1}{c|}{610.327} & \multicolumn{1}{c|}{80} & 81.80 & \multicolumn{1}{c|}{1144.32} & \multicolumn{1}{c|}{152} & 83.31 & 1171.81 & 135 \\
		\hline
		4\_8\_50 & 102.25 & \multicolumn{1}{c|}{559.198} & \multicolumn{1}{c|}{129} & 101.79 & \multicolumn{1}{c|}{598.106} & \multicolumn{1}{c|}{95} & 101.41 & 5070.40 & 122 \\
		\hline
		5\_8\_50 & 82.71 & \multicolumn{1}{c|}{601.796} & \multicolumn{1}{c|}{85} & 83.22 & \multicolumn{1}{c|}{703.881} & \multicolumn{1}{c|}{107} & 82.90 & 1748.01 & 199 \\
		\hline
		6\_8\_50 & 81.55 & \multicolumn{1}{c|}{703.951} & \multicolumn{1}{c|}{126} & 81.66 & \multicolumn{1}{c|}{789.632} & \multicolumn{1}{c|}{116} & 81.49 & 1004.17 & 120 \\
		\hline
		7\_8\_50 & 93.24 & \multicolumn{1}{c|}{555.179} & \multicolumn{1}{c|}{90} & 93.29 & \multicolumn{1}{c|}{666.236} & \multicolumn{1}{c|}{80} & 93.20 & 733.42 & 113 \\
		\hline
		8\_8\_50 & 83.49 & \multicolumn{1}{c|}{802.583} & \multicolumn{1}{c|}{117} & 81.43 & \multicolumn{1}{c|}{820.124} & \multicolumn{1}{c|}{133} & 80.40 & 803.35 & 118 \\
		\hline
		9\_8\_50 & 82.81 & \multicolumn{1}{c|}{884.723} & \multicolumn{1}{c|}{134} & 80.26 & \multicolumn{1}{c|}{1040.52} & \multicolumn{1}{c|}{142} & 80.39 & 1144.87 & 127 \\
		\hline
		10\_8\_50 & 80.94 & \multicolumn{1}{c|}{830.715} & \multicolumn{1}{c|}{130} & 80.27 & \multicolumn{1}{c|}{990.318} & \multicolumn{1}{c|}{138} & 80.01 & 1329.35 & 164 \\
		\hline
		\textit{\textbf{Average Values}} & 87.36 & 655.47 & 107.00 & 86.87 & 793.25 & 116.20 & 86.70 & 1469.18 & 136.80 \\
		\hline
	\end{tabular}
	\label{tab:iterlimit}%
\end{table}%

The values of the objective function and computation times are given in Table \ref{tab:iterlimit}, where we observe that the average number of iterations before convergence increases with the higher values of $iterlimit$, resulting with higher average run times. When $iterlimit$ is set to lower values, the algorithm starts to use higher values of $\rho$ in earlier iterations and the quality of the solutions worsen. Since we would not like to sacrifice solution quality and computational times are acceptable at this level of the parameter, we fix $iterlimit$ at 100 in our experiments.

\end{document}